\RequirePackage{fix-cm}
\documentclass[twoside,british,times,3p,fleqn]{elsarticle}
\usepackage[T1]{fontenc}
\usepackage[latin9]{inputenc}
\usepackage{babel}
\usepackage{array}
\usepackage{bm}
\usepackage{multirow}
\usepackage{amsmath}
\usepackage{amsthm}
\usepackage{amssymb}
\usepackage{mathdots}
\usepackage{graphicx}
\usepackage{esint}
\usepackage[unicode=true,
 bookmarks=true,bookmarksnumbered=false,bookmarksopen=false,
 breaklinks=false,pdfborder={0 0 1},backref=false,colorlinks=false]
 {hyperref}
\hypersetup{pdftitle={On Generalisation of Isotropic Central Difference for Higher Order Approximation of Fractional Laplacian},
 pdfauthor={P.H. Lam},
 pdfkeywords={Finite Difference, Fractional Laplacian, Fractional Derivative, Higher Order Approximation, Smolyak}}

\makeatletter

\providecommand{\tabularnewline}{\\}

\numberwithin{equation}{section}
\numberwithin{figure}{section}
\numberwithin{table}{section}

\journal{Journal of Computational Physics}

\makeatother

\begin{document}

\begin{frontmatter}{}

\title{On Generalisation of Isotropic Central Difference for Higher Order
Approximation of Fractional Laplacian}

\author{P.H.~Lam\corref{cor1}}

\ead{puiholam2-c@my.cityu.edu.hk}

\author{H.C.~So}

\ead{hcso@ee.cityu.edu.hk}

\cortext[cor1]{Corresponding author}

\address{Department of Electrical Engineering, City University of Hong Kong,
Kowloon, Hong Kong SAR}
\begin{abstract}
The study of generalising the central difference for integer order
Laplacian to fractional order is discussed in this paper. Analysis
shows that, in contrary to the conclusion of a previous study, difference
stencils evaluated through fast Fourier transform prevents the convergence
of the solution of fractional Laplacian. We propose a composite quadrature
rule in order to efficiently evaluate the stencil coefficients with
the required convergence rate in order to guarantee convergence of
the solution. Furthermore, we propose the use of generalised higher
order lattice Boltzmann method to generate stencils which can approximate
fractional Laplacian with higher order convergence speed and error
isotropy. We also review the formulation of the lattice Boltzmann
method and discuss the explicit sparse solution formulated using Smolyak's
algorithm, as well as the method for the evaluation of the Hermite
polynomials for efficient generation of the higher order stencils.
Numerical experiments are carried out to verify the error analysis
and formulations.
\end{abstract}

\end{frontmatter}{}

\section{Introduction}

Fractional Laplacian has found important applications in modelling
physics, which include low frequency approximation of the Zener fractional
wave model \citep{treeby2010modeling,holm2014comparison}, fractional
diffusion model for porous medium \citep{depablo2011afractional},
diffusion reaction \citep{volpert2013frontsin}, and more can be found
in the review \citep{lischke2020whatis}. Therefore, numerical methods
for fractional Laplacian is also a hot topic.

There are many definitions of fractional Laplacian. In particular,
one that is defined for non-local boundary condition needs that the
function undertaking the Laplacian to be defined in everywhere including
the exterior of the domain under consideration. Another definition
based on the eigenvalue of a boundary problem only requires the function
to be defined within the boundary. For a full review of various definitions
of the fractional Laplacian and numerical solutions, refer to \citep{lischke2020whatis}.
In this paper, we consider the Riesz type, which belongs to the former
definition, where the function must be defined in the entirety of
real number space. We follow a similar approach as the previous work
\citep{hao2021fractional}, where we seek a filter with a Fourier
spectrum that closely resembles that of Riesz fractional Laplacian.
Surprisingly, this approach has not yet been addressed in the literature
\citep{lischke2020whatis} even though it dates as back to 2005 \citep{ilic2005numerical},
and perhaps much earlier.

In the previous work \citep{hao2021fractional}, it was concluded
that evaluation of the stencil coefficients by fast Fourier transform
(FFT) is sufficient to achieve second and higher convergence. However,
we will show that this is not the case because the spectrum of the
stencil does not converge properly, resulting in a constant error,
regardless of the nodal spacing. Moreover, the numerical examples
only include low number of nodes, restricting the approximation error
of the solution to a much higher order. This hides the fact that the
spectral property of the stencil does not converge properly. A more
sophisticated approach is called for to solve for those stencil coefficients
in order to guarantee convergence order of the solution. Other than
addressing this issue, we have made a number of improvements to it
by incorporating the generalised lattice Boltzmann method. However,
we do not discuss the application of this stencil in dynamic problems
as it has already been addressed in \citep{hao2021fractional} and
many other papers.

The paper is organised as follows. First, we define notations vital
for understanding of the mathematical expressions in this paper in
Subsection \ref{subsec:Notations}. In Section \ref{sec:Central-Difference-based},
we review the generalised lattice Boltzmann method in details with
additional work of our own. In Section \ref{sec:Generalisation-to-Fractional},
we discuss how the lattice Boltzmann method can be generalised to
fractional order for approximating the fractional Laplacian. Error
analysis and convergence guarantee of various methods incorporated
are presented. In Section \ref{sec:Numerical-Experiment}, we demonstrate
that, in practice, the use of FFT for obtaining stencils may prevent
the convergence of the solution and that our solution successfully
averts this issue. And finally in the conclusion, we remark on the
applied methods and suggest future study.

We would also like to highlight the contributions of this paper as
follows:
\begin{itemize}
\item Solution and rule set for automatic generation of central difference
with arbitrary order of convergence and error isotropy
\item Error analysis of the generalised lattice Boltzmann is presented using
the multinomial formula for understanding of the relationship between
the convergence and error isotropy order with the order of Hermite
quadrature
\item Explicit solutions, in terms of Lagrange polynomials, for sparse multidimensional
lattice Boltzmann stencils using modified Smolyak method
\item More efficient method of computing the Hermite polynomials required
for higher order stencil compared to \citep{mattila2014highaccuracy}
\item Convergence analysis of the tanh-sinh rule for inverse discrete-time
Fourier transform (DTFT)
\item Convergence analysis of generalised Filon method to ensure that the
fractional central difference converges at a rate required for the
convergence of the solution at expected rate
\end{itemize}

\subsection{\label{subsec:Notations}Notations}

Here, we define the commonly used notations throughout the paper.
Tensors are denoted by bold italic characters, and an element of an
$n$ rank tensor is denoted by a tensor with a superscript of a vector
consisting of $n$ elements. For example, indexing of the tensor $\boldsymbol{T}$
with the vector indices $\mathbf{a}$ is $\boldsymbol{T}^{\mathbf{a}}$.
The number of dimensions is denoted by $d$. We use the probabilist's
definition of Hermite polynomials given by 
\begin{equation}
H_{n}=\left(-1\right)^{n}\mathrm{e}^{\frac{x^{2}}{2}}\frac{\mathrm{d}^{n}}{\mathrm{d}x^{n}}\mathrm{e}^{-\frac{x^{2}}{2}}.
\end{equation}
The $n$ rank Hermite tensor of dimension $d^{n}$ in the Cartesian
coordinate system is defined as
\begin{equation}
\boldsymbol{H}_{n}^{\mathbf{a}}\left(\mathbf{x}\right)=\prod_{k=0}^{d-1}H_{s\left(\mathbf{a}|k\right)}\left(\mathbf{x}_{k}\right),\label{eq:hermite tensor}
\end{equation}
where $\mathbf{a}$ is a vector which consists of $n$ integer indices
for $0\le\mathbf{a}_{j}<d$, and the parity function $s\left(\cdot\right)$
is defined as
\begin{equation}
s\left(\mathbf{a}|k\right)=\sum_{j=0}^{n-1}\delta_{\mathbf{a}_{j},k}.
\end{equation}
We also use the vector subscript on a function to define the order
of derivative applied for each respective dimension of the vector
\begin{equation}
\left\{ \frac{\partial^{\otimes n}}{\partial\mathbf{x}^{\otimes n}}f\left(\mathbf{x}\right)\right\} ^{\mathbf{a}}=f_{\mathbf{s}\left(\mathbf{a}\right)}\left(\mathbf{x}\right),
\end{equation}
where $\left\{ \mathbf{s}\left(\mathbf{a}\right)\right\} _{k}=s\left(\mathbf{a}|k\right)$.
The Laplacian operator is defined by the symbol $\Delta$, and the
higher order Laplacian is defined as
\begin{equation}
\Delta^{m}f\left(\mathbf{x}\right)=\left(\sum_{k=0}^{d-1}\frac{\partial^{2}}{\partial\mathbf{x}_{k}^{2}}\right)^{m}f\left(\mathbf{x}\right).
\end{equation}
Moreover, we define the fractional Laplacian of $f$ as
\begin{equation}
-\left(-\Delta\right)^{\frac{\alpha}{2}}f\left(\mathbf{x}\right)=-\mathcal{F}^{-1}\left\{ \left|\mathbf{k}\right|^{\alpha}\mathcal{F}\left\{ f\right\} \left(\mathbf{k}\right)\right\} \left(\mathbf{x}\right),
\end{equation}
where $\mathcal{F}$ is the Fourier transform operator.

\section{Central Difference based on Hermite Polynomials with Higher Order
Isotropy\label{sec:Central-Difference-based}}

First we review the derivation of lattice Boltzmann method in the
form of a central difference operator to approximate an arbitrary
order of derivative. The difference operator $\boldsymbol{D}_{n}$
of order $n$ in the lattice Boltzmann method is defined by the $d$
rank tensor contraction \citep{philippi2006fromthe}
\begin{equation}
f_{\mathbf{s}\left(\mathbf{a}\right)}\left(\mathbf{x}\right)\sim\boldsymbol{D}_{n}^{\mathbf{a}}\cdot\boldsymbol{F}=\left(\frac{a}{h}\right)^{n}\sum_{\substack{\mathbf{i}_{j}\in\left\{ -N_{h},\ldots,N_{\mathrm{h}}\right\} ,\\
j=0,\ldots,d-1
}
}w_{\mathbf{i}}^{N_{q}}\boldsymbol{H}_{n}^{\mathbf{a}}\left(\mathbf{v}_{\mathbf{i}}\right)f\left(\mathbf{x}+\mathbf{x}_{\mathbf{i}}\right),\label{eq:latticeboltzman def}
\end{equation}
where $\mathbf{x}_{\mathbf{i}}=\mathbf{i}h$ are the grid points with
equal spatial distance $h$, $\mathbf{v}_{\mathbf{i}}=a\mathbf{i}$
are the grid points scaled by a constant factor $a>0$, and $w_{\mathbf{i}}^{N_{q}}$
are quadrature weights at the grid points. These weights are invariant
in the sense that they are the same when swapping indices, and they
satisfy the norm preservation of the continuous counterpart of inner
products of Hermite polynomials, for all $n\le N_{q}$ and unique
combinations of $\mathbf{a}$ such that $\sum_{k}\mathbf{a}_{k}=n$
with each $\mathbf{a}_{k}\ge\mathbf{a}_{k+1}$, as such
\begin{equation}
\sum_{\mathbf{i}\in\bm{\Omega}}w_{\mathbf{i}}^{N_{q}}\left(\boldsymbol{H}_{n}^{\mathbf{a}}\left(\mathbf{v}_{\mathbf{i}}\right)\right)^{2}=\int_{\mathrm{R}^{d}}\prod_{k=0}^{d-1}H_{s\left(\mathbf{a}|k\right)}^{2}\left(\mathbf{x}_{k}\right)\omega\left(\mathbf{x}_{k}\right)\,\mathrm{d}\mathbf{x},\label{eq:quadrature problem}
\end{equation}
where $\Omega$ is the space of indices for non-zero weights to simplify
the multi-index summation, and $\omega\left(x\right)=\frac{1}{\sqrt{2\pi}}\mathrm{e}^{-\frac{x^{2}}{2}}$
is the Hermite weight function. The purpose of the scaling factor
$a$, and the number of of weights per dimension $N_{h}$ required
will be discussed later in Subsection \ref{subsec:1D-Quadrature-Problem}.

If the weights satisfy these norm preservation conditions, they also
satisfy the orthogonality conditions for inner products of Hermite
polynomials. To understand how, first, consider the 1D weights, due
to symmetry of the weights, when the inner product consists of polynomials
of different parities, the integrand is odd, and they cancel each
other from the summation of positive and negative sides. When the
parities are the same but the polynomials are not of the same order,
the integrand is a polynomial consisting of only monomials of even
degrees. These monomials can be written as a linear combination of
each of the squared polynomials because one can show that the lower
triangle matrix formed from writing the squared polynomials as linear
combinations of the monomials has a non-zero determinant. In fact,
the inverse is given explicitly given by (14) in \citep{qi2021someproperties}.
Since each integral of each of these squared polynomials is exactly
equal to the discrete sum, the integral of the linear combinations
will also yield the same result. Therefore, the discrete sum also
satisfies the orthogonality of the continuous counterpart. In the
multidimensional case, because the integrand is simply a product of
the 1D squared polynomials in each dimension, satisfying all combinations
of which the total order is lower or equal to $N_{q}$ is enough.
Moreover, since the weights are invariant, permutations of the orders
of the polynomial in each dimension lead to the same sum, that is,
the sum is also invariant. Therefore, we only need to satisfy unique
combinations of degree numbers which sum to $\le N_{q}$, regardless
of the permutation.

To understand the motivation of this approximation, first, we suppose
$N_{q}\rightarrow\infty$ so that \eqref{eq:latticeboltzman def}
can be rewritten as
\begin{align}
\boldsymbol{D}_{n}^{\mathbf{a}}\cdot\boldsymbol{F}=\left(\frac{a}{h}\right)^{n} & \int_{\mathrm{R}^{d}}\boldsymbol{H}_{n}^{\mathbf{a}}\left(\mathbf{v}\right)\omega\left(\mathbf{v}\right)f\left(\mathbf{x}-\frac{h}{a}\mathbf{v}\right)\,\mathrm{d}\mathbf{v}\\
=\left(-\frac{a}{h}\right)^{n} & \int_{\mathrm{R}^{d}}\prod_{k=0}^{d-1}\left(\frac{\mathrm{d}^{s\left(\mathbf{a}|k\right)}}{\mathrm{d}\mathbf{v}_{k}^{s\left(\mathbf{a}|k\right)}}\mathrm{e}^{-\frac{\mathbf{v}_{k}^{2}}{2}}\right)f\left(\mathbf{x}-\frac{h}{a}\mathbf{v}\right)\,\mathrm{d}\mathbf{v}\\
= & \int_{\mathrm{R}^{d}}\mathrm{e}^{-\frac{\mathbf{v}_{k}^{2}}{2}}\left.\left\{ \frac{\partial^{\otimes n}}{\partial\bm{\xi}^{\otimes n}}f\left(\bm{\xi}\right)\right\} _{\mathbf{a}}\right|_{\bm{\xi}=\mathbf{x}-\frac{h}{a}\mathbf{v}}\,\mathrm{d}\mathbf{v}\\
= & \mathcal{F}^{-1}\left\{ \mathrm{e}^{-\frac{1}{2}\left|\frac{h}{a}\mathbf{k}\right|^{2}}\prod_{k=0}^{d-1}\left(\mathrm{i}\mathbf{k}_{k}\right)^{s\left(\mathbf{a}|k\right)}\mathcal{F}\left\{ f\right\} \left(\mathbf{k}\right)\right\} \left(\mathbf{x}\right),
\end{align}
where we have assumed that $f$ converges everywhere when expanded
as a Taylor's series about any point $\mathbf{x}$ and that $f$ vanishes
in the infinity when multiplied by the Gaussian function. This shows
that the approximation is the derivative of $f$ spectrally low-pass
filtered by a Gaussian function. Expanding the Gaussian function about
$0$ leads to
\begin{equation}
\boldsymbol{D}_{n}^{\mathbf{a}}\cdot\boldsymbol{F}=\sum_{k=0}^{\infty}\frac{1}{2^{k}k!}\left(\frac{h}{a}\right)^{2k}\Delta^{k}f_{\mathbf{s}\left(\mathbf{a}\right)}\left(\mathbf{x}\right),\label{eq:lattice boltzman infinite order error analysis}
\end{equation}
which shows that the operator gives the second order approximation
of the derivative, with infinite order of error isotropy.

We can show that a finite $N_{q}$ truncates this series to $k=N_{q}-n$,
provided that $N_{q}\ge n$, with higher order terms being unknown,
and thus the approximation can still be second order accurate, and
the error can still be isotropic up to order $2\left(N_{q}-n\right)$.
Substituting the Taylor's series expansion of $f$ about $\mathbf{x}$,
expressed in the multinomial form, into \eqref{eq:latticeboltzman def}
gives
\begin{align}
\boldsymbol{D}_{n}^{\mathbf{a}}\cdot\boldsymbol{F} & =\left(\frac{a}{h}\right)^{n}\sum_{m=0}^{2N_{q}-n}\frac{1}{m!}\int_{\mathrm{R}^{d}}\left(\prod_{k=0}^{d-1}H_{s\left(\mathbf{a}|k\right)}\left(\mathbf{v}_{k}\right)\right)\left(\sum_{k=0}^{d-1}\frac{h}{a}\mathbf{v}_{k}\frac{\partial}{\partial\mathbf{x}_{k}}\right)^{m}f\left(\mathbf{x}\right)\,\mathrm{d}\mathbf{v}+O\left(h^{2N_{p}-n+1}\right)\label{eq:lattice boltzman finite error order}\\
 & \sim\left(\frac{a}{h}\right)^{n}\sum_{m=0}^{2N_{q}-n}\left(\frac{h}{a}\right)^{m}\sum_{\sum_{k=0}^{d-1}\mathbf{b}_{k}=m}\frac{f_{\mathbf{b}}\left(\mathbf{x}\right)}{\left(2\pi\right)^{\frac{d}{2}}}\prod_{k=0}^{d-1}\int_{-\infty}^{\infty}\frac{\mathbf{v}_{k}^{\mathbf{b}_{k}}}{\mathbf{b}_{k}!}H_{s\left(\mathbf{a}|k\right)}\left(\mathbf{v}_{k}\right)\omega\left(\mathbf{v}_{k}\right)\,\mathrm{d}\mathbf{v}_{k},
\end{align}
where the summation limits refers to all possible combinations of
$\mathbf{b}$ which satisfy the condition. The integral is only non-zero
when $\mathbf{b}_{k}-s\left(\mathbf{a}|k\right)$, and thus $m-n$,
is even and non-negative. This leaves us with only the combinations
of the squared terms. Moreover, the integral with respect to each
dimension results in $2^{-l}/l!$, where $l=\left(\mathbf{b}_{k}-s\left(\mathbf{a}|k\right)\right)/2$.
Therefore, we can rewrite the difference as
\begin{align}
\boldsymbol{D}_{n}^{\mathbf{a}}\cdot\boldsymbol{F} & =\sum_{m=0}^{N_{q}-n}\frac{1}{2^{m}m!}\left(\frac{h}{a}\right)^{2m}\sum_{\sum_{k=0}^{d-1}\mathbf{b}_{k}=m}\prod_{k=0}^{d-1}\frac{m!}{\mathbf{b}_{k}!}\frac{\partial^{2\mathbf{b}_{k}}}{\partial\mathbf{x}_{k}^{2\mathbf{b}_{k}}}f_{\mathbf{s}\left(\mathbf{a}\right)}\left(\mathbf{x}\right)+O\left(h^{2\left(N_{q}-n+1\right)}\right)\\
 & =\left(1+\sum_{m=1}^{N_{q}-n}\frac{1}{2^{m}m!}\left(\frac{h}{a}\right)^{2m}\left(\sum_{k=0}^{d-1}\frac{\partial^{2}}{\partial\mathbf{x}_{k}^{2}}\right)^{m}\right)f_{\mathbf{s}\left(\mathbf{a}\right)}\left(\mathbf{x}\right)+O\left(h^{2\left(N_{q}-n+1\right)}\right),\label{eq:lattice boltzman final error}
\end{align}
which is basically \eqref{eq:lattice boltzman infinite order error analysis}
but with terms higher order than $N_{q}-n$ being unknowns. This shows
that as long as $N_{q}\ge n$, the error is second order. However,
it is only isotropic up to order $2\left(N_{q}-n\right)$, for $N_{q}>n$.

\subsection{\label{subsec:Higher-Order-Approximation}Higher Order Approximation}

Because of the unique properties of the Taylor's series coefficients
of the Gaussian function, the higher order accurate difference stencils
can in fact be explicitly expressed in terms of the quadrature weights
and higher order Hermite polynomials. Here we show that an arbitrary
order $2N_{c}$ accurate difference method is given by
\begin{align}
\boldsymbol{D}_{n}^{\mathbf{a}}\cdot\boldsymbol{F} & =\left(\frac{a}{h}\right)^{n}\sum_{\mathbf{i}\in\bm{\Omega}}w_{\mathbf{i}}^{N_{q}}\sum_{j=0}^{N_{c}-1}\frac{\left(-1\right)^{j}}{2^{j}j!}\boldsymbol{H}_{n,j}^{\mathbf{a}}\left(\mathbf{v}_{\mathbf{i}}\right),\label{lattice boltzman higher order}\\
 & =f_{\mathbf{s}\left(\mathbf{a}\right)}\left(\mathbf{x}\right)+O\left(h^{2N_{c}}\right),
\end{align}
where $N_{c}$ satisfies $N_{q}\ge n+2\left(N_{c}-1\right)$, and
we define the Laplacian Hermite polynomials as
\begin{equation}
\boldsymbol{H}_{n,j}^{\mathbf{a}}\left(\mathbf{x}\right)=\left(-1\right)^{n}\mathrm{e}^{\sum_{k=0}^{d-1}\frac{\mathbf{x}_{k}^{2}}{2}}\Delta^{j}\prod_{k=0}^{d-1}\frac{\mathrm{d}^{\mathbf{a}_{k}}}{\mathrm{d}\mathbf{x}_{k}^{\mathbf{a}_{k}}}\mathrm{e}^{-\frac{\mathbf{x}_{k}^{2}}{2}}.\label{eq:hermite laplacian of any order}
\end{equation}
From \eqref{eq:lattice boltzman final error}, we deduce that \eqref{lattice boltzman higher order}
leads to
\begin{multline}
\boldsymbol{D}_{n}^{\mathbf{a}}\cdot\boldsymbol{F}=\sum_{k=0}^{N_{c}-1}\left(\sum_{l=0}^{k}\frac{1}{2^{l}l!}\frac{\left(-1\right)^{k-l}}{2^{k-l}\left(k-l\right)!}\right)\left(\frac{h}{a}\right)^{2k}\Delta^{k}f_{\mathbf{s}\left(\mathbf{a}\right)}\left(\mathbf{x}\right)\\
+\sum_{k=N_{c}}^{N_{q}-n-N_{c}+1}\left(\sum_{l=0}^{N_{c}-1}\frac{\left(-1\right)^{l}}{2^{l}l!}\frac{1}{2^{k-l}\left(k-l\right)!}\right)\left(\frac{h}{a}\right)^{2k}\Delta^{k}f_{\mathbf{s}\left(\mathbf{a}\right)}\left(\mathbf{x}\right)+O\left(h^{2\left(N_{q}-n-N_{c}+2\right)}\right).\label{eq:detailed error terms}
\end{multline}
The coefficient for $k<N_{c}$ can be simplified as
\begin{align}
\sum_{l=0}^{k}\frac{1}{2^{l}l!}\frac{\left(-1\right)^{k-l}}{2^{k-l}\left(k-l\right)!} & =\frac{\left(-1\right)^{k}}{2^{k}k!}\sum_{l=0}^{k}\frac{\left(-1\right)^{l}k!}{l!\left(k-l\right)!}\\
 & =\frac{\left(-1\right)^{k}}{2^{k}k!}\sum_{l=0}^{k}\binom{k}{l}\left(-1\right)^{l}\left.x^{k-l}\right|_{x=1}\\
 & =\frac{\left(-1\right)^{k}}{2^{k}k!}\left.\left(x-1\right)^{k}\right|_{x=1}=0.
\end{align}
This proves that the error terms of order lower than $2N_{c}$ are
cancelled out, provided that $N_{q}$ is sufficiently large. The order
of isotropy of the error follows from before, and it requires $N_{q}$
to satisfy $N_{q}\ge2N_{c}+n-1$ for at least an error term to be
isotropic.

\subsection{\label{subsec:1D-Quadrature-Problem}1D Grid Quadrature Problem for
Hermite Weight Function}

Recall that the lattice Boltzmann method has turned the derivative
problem into an integral problem where we need to find coefficients
such that \eqref{eq:quadrature problem} is satisfied. In the literature,
these weights are commonly solved in the implicit matrix form \citep{philippi2006fromthe}.
We shall review the formulation of this problem in 1D and higher dimensions.
Moreover, we show that the solution can be explicitly expressed in
terms of Lagrange polynomials.

In 1D, this problem is straightforward to formulate as we can simply
choose consecutive nodes on the grid so that the degree of freedom
matches the number of conditions to satisfy. Because the squared Hermite
polynomials are symmetric, the weights are also symmetric by Sobolev's
invariant theorem \citep{krommer1998computational}. Therefore, the
weights on the negative side, in fact, do not contribute to the degree
of freedom. Then, in matrix form, the linear system of equations consisting
of \eqref{eq:quadrature problem} for various $n$ is written as
\begin{equation}
\begin{bmatrix}H_{0}^{2}\left(0\right) & H_{0}^{2}\left(a\right) & \cdots & H_{0}^{2}\left(N_{q}a\right)\\
H_{1}^{2}\left(0\right) & H_{1}^{2}\left(a\right) & \cdots & H_{1}^{2}\left(N_{q}a\right)\\
\vdots & \vdots & \ddots & \vdots\\
H_{N_{q}}^{2}\left(0\right) & H_{N_{q}}^{2}\left(a\right) & \cdots & H_{N_{q}}^{2}\left(N_{q}a\right)
\end{bmatrix}\begin{bmatrix}1 & 0 & \cdots & 0\\
0 & 2 & \iddots & 0\\
\vdots & \iddots & \ddots & \vdots\\
0 & \cdots & 0 & 2
\end{bmatrix}\begin{bmatrix}w_{0}\\
w_{1}\\
\vdots\\
w_{N_{q}}
\end{bmatrix}=\begin{bmatrix}0!\\
1!\\
\vdots\\
N_{q}!
\end{bmatrix}.\label{eq:1d hermite grid quadrature linear system}
\end{equation}
Alternatively, one can also express this in terms of monomials and
the Hermite moments as
\begin{equation}
\begin{bmatrix}1 & 1 & \cdots & 1\\
0 & a^{2} & \cdots & \left(N_{q}a\right)^{2}\\
\vdots & \vdots & \ddots & \vdots\\
0 & \left(a\right)^{2N_{q}} & \cdots & \left(N_{q}a\right)^{2N_{q}}
\end{bmatrix}\begin{bmatrix}1 & 0 & \cdots & 0\\
0 & 2 & \iddots & 0\\
\vdots & \iddots & \ddots & \vdots\\
0 & \cdots & 0 & 2
\end{bmatrix}\begin{bmatrix}w_{0}\\
w_{1}\\
\vdots\\
w_{N_{q}}
\end{bmatrix}=\begin{bmatrix}1\\
\frac{2!}{2\left(1!\right)}\\
\vdots\\
\frac{\left(2N_{q}\right)!}{2^{N_{q}}N_{q}!}
\end{bmatrix}.
\end{equation}
The right hand side is the so-called double factorial. Solving this
formulation is less stable but save on costs of computing the polynomials.
With the extra variable $a$, one of the weights can be eliminated.
Ideally, one would want the weight eliminated to be the furthest weight
from the centre but this is not always possible. Numerical solutions
suggest that $w_{N_{q}}$ can only be eliminated for even $N_{q}$.
Therefore, for the rest of the paper, we follow this conjecture. However,
solving for $a$ in the implicit form is although possible, it is
likely expensive. Here, we show that it can be reduced into a root
finding problem by explicitly expressing the weights in terms of the
Hermite polynomials.

To formulate the explicit solution, note that the squared Hermite
polynomials can be interpolated exactly by $2N_{q}+1$ points Lagrange
polynomials. The integral itself can be evaluated exactly using Gauss
quadrature with $N_{q}+1$ points. Therefore, the right hand side
of \eqref{eq:quadrature problem} can be rewritten as
\begin{equation}
\sum_{k=-N_{q}}^{N_{q}}H_{n}^{2}\left(ak\right)\sum_{i=0}^{N_{q}}q_{i}l_{k}^{N_{q}}\left(\frac{\xi_{i}}{a}\right),\label{eq:rhs 1d lagrange}
\end{equation}
where $q_{i}$ and $\xi_{i}$ are the Gauss-Hermite quadrature weights
and abscissae respectively, and each $l_{k}^{n}$ is the $2n+1$ symmetric
grid point Lagrange polynomial basis defined as
\begin{equation}
l_{k}^{n}\left(x\right)=\prod_{\substack{j=-n\\
j\neq k
}
}^{n}\frac{x-j}{k-j}.
\end{equation}
Comparing the left and right hand sides, we see that
\begin{equation}
w_{k}=\sum_{i=0}^{N_{q}}q_{i}l_{k}^{N_{q}}\left(\frac{\xi_{i}}{a}\right).\label{eq:hermite grid quadrature 1D solution}
\end{equation}
The Gaussian weights and abscissa are most easily obtained by solving
the eigenvalue problem for the tridiagonal matrix of dimension $\mathbb{R}^{N_{q}\times N_{q}}$
given by \citep{golub1969calculation}
\begin{equation}
\begin{bmatrix}0 & 1\\
1 & 0 & \sqrt{2}\\
 & \sqrt{2} & \ddots & \ddots\\
 &  & \ddots & \ddots & \sqrt{N_{q}-1}\\
 &  &  & \sqrt{N_{q}-1} & 0
\end{bmatrix}.
\end{equation}
The diagonal elements are obtained from the recurrence relationship
of the Hermite polynomials. The eigenvalues are the roots of $H_{N_{q}+1}$
and thus Gauss quadrature abscissae, and the weights are given by
$v_{k}^{2}$, where $v_{k}$ is the first element of the eigenvector
corresponding to the $k$-th eigenvalue. Again, Sobolev's theorem
ensures that these weights and abscissae are symmetric, so that $q_{i}=q_{N_{q}-i}$,
and $\xi_{i}=\xi_{N_{q}-i}$. Moreover, since $l_{N_{q}}\left(x\right)=l_{-N_{q}}\left(-x\right)$,
one can rewrite the linear system of \eqref{eq:rhs 1d lagrange} in
matrix form, similarly as \eqref{eq:1d hermite grid quadrature linear system}.
The cancellation of the matrix consisting of the Hermite polynomials
on both sides confirms that \eqref{eq:hermite grid quadrature 1D solution}
is indeed the solution of \eqref{eq:1d hermite grid quadrature linear system}.
This also validates the requirement of $2N_{q}+1$ symmetric nodes
even though there are only $N_{q}+1$ conditions to satisfy, which
makes it appears to require only $N_{q}+1$ weights.

With the closed-form solution, the scaling constant $a$ can be solved
by setting $w_{k}^{N_{q}}=0$, resulting in the root finding problem
of solving for $a$ such that
\begin{equation}
\sum_{i=0}^{N_{q}}q_{i}\xi_{i}\left(\xi_{i}+ak\right)\prod_{\substack{j=1\\
j\neq k
}
}^{N_{q}}\left(\frac{\xi_{i}^{2}}{j^{2}}-a^{2}\right)=0.
\end{equation}
Using the symmetry of the quadrature weights and abscissae and $\xi_{N_{q}/2}=0$
for even $N_{q}$, the problem can be reduced to
\begin{align}
\sum_{i=0}^{\left\lceil \frac{N_{q}}{2}\right\rceil -1}q_{i}\prod_{\substack{j=0\\
j\neq k
}
}^{N_{q}}\left(\xi_{i}^{2}-a^{2}j^{2}\right) & =0.\label{eq:weight elimination condition}
\end{align}
Alternatively, one can obtain the polynomial monomial integer coefficients
first by convolution and evaluating the moments, resulting in
\begin{equation}
\sum_{i=1}^{N_{q}}c_{i}\left(a^{2}\right)^{i-1}=0,
\end{equation}
where $c_{i}=\frac{\left(2i\right)!}{2^{i}i!}\mathcal{F}^{-1}\left\{ \prod_{j=1,j\neq k}^{N_{q}}\left(1-j^{2}\mathrm{e}^{-\mathrm{i}\omega}\right)\right\} \left(i-1\right).$
It is not exactly known whether there is a positive real solution
for $a$ when $k=N_{q}$, other than the aforementioned conjecture
when $N_{q}$ is even, but one can show that when $w_{n}^{n}=0$,
$w_{m}^{n-1}=w_{m}^{n}$ for $m=0,\,\ldots,\,n-1$. In fact, one can
also show that all those conditions for various $m$ are equivalent
to \eqref{eq:weight elimination condition} when $k=n$. Therefore,
\eqref{eq:weight elimination condition} is the only condition we
have to solve when seeking an appropriate $a$ for the elimination
of the last coefficient.

\subsection{Multidimensional Quadrature Problem\label{subsec:Multidimensional-Quadrature-Prob}}

In the multidimensional case, one can easily show that the tensor
product of quadrature weights satisfy \eqref{eq:quadrature problem}
since each dimension can be separated by multiplication. However,
from \eqref{eq:lattice boltzman finite error order}, it is understood
that the quadrature rule only needs to satisfy the exactness of the
integral over a homogeneous polynomial of degree $N_{q}$. That is
the maximum of the sum of monomial orders of all dimensions only equals
$N_{q}$, but not $dN_{q}$. Although the tensor product of quadrature
weights satisfies all the terms in a monomial tensor, it is dense
and contains many more coefficients than necessary. This is commonly
referred to as curse of dimensionality, and it can accumulate errors
from the additional summations. To solve this problem, the Smolyak's
method can be applied to produce sparse weights which satisfy the
orthogonality requirement stated by \eqref{eq:quadrature problem}
\citep{gerstner1998numerical}. Before we discuss the explicit solution
obtained from Smolyak's method, let us discuss the implicit solution.

From \eqref{eq:quadrature problem}, the number of conditions is given
by the triangle number $\left(\left(2+\left\lfloor \frac{N_{q}}{2}\right\rfloor \right)\left(1+\left\lfloor \frac{N_{q}}{2}\right\rfloor \right)+\left(\left\lceil \frac{N_{q}}{2}\right\rceil +1\right)\left\lceil \frac{N_{q}}{2}\right\rceil \right)/2$
in 2D. The triangle number requires that the grid spans at least $N_{q}$
number of nodes on each side, the same as the 1D case. Recall that
each coefficient $w_{i,j}$ equals $w_{j,i}$ because of the invariance
theorem. Thus, we choose only nodes from a right triangle formed by
$\frac{N_{q}^{2}}{2}$ positive nodes. Further halving this right
triangle as an isosceles triangle results in this triangle number
of nodes. The nodes on the axes and the diagonal nodes only repeat
4 times while the nodes in between are repeated 8 times. Therefore,
prioritizing nodes on the axes and diagonals results in fewer nodes.
This requires shifting the nodes in the other half of the isosceles
triangle towards the diagonal axis. For the nodes to be more evenly
spread out, we adopt the following rule: $w_{i,j}$ are non-zero for
$i=0,\,\ldots,\,\left\lfloor \frac{N_{q}}{2}\right\rfloor $, $j=0,\,\ldots,\,i$,
and for $i=\left\lfloor \frac{N_{q}}{2}\right\rfloor +1,\,\ldots N_{q}$,
$j=0,\,\left\lfloor \frac{i}{N_{q}-i}\right\rfloor ,\,\left\lfloor \frac{2i}{N_{q}-i}\right\rfloor ,\,\ldots$.
This gives us the most optimised scheme for the 2D case.

For example, for the case $N_{q}=4$, we have the following set of
nodes $\mathbf{i}\in\left\{ \left(0,0\right),\left(1,0\right),\left(1,1\right),\left(2,0\right),\left(2,1\right),\left(2,2\right),\right.$

\noindent $\left.\left(3,0\right),\left(3,3\right),\left(4,0\right)\right\} .$
The orthogonality conditions to satisfy are then
\begin{equation}
\begin{bmatrix}\left(H_{4}\left(0\right)H_{0}\left(0\right)\right)^{2} & \left(H_{4}\left(a\right)H_{0}\left(0\right)\right)^{2} & \cdots & \left(H_{4}\left(3a\right)H_{0}\left(3a\right)\right)^{2} & \left(H_{4}\left(2a\right)H_{0}\left(a\right)\right)^{2}\\
\left(H_{3}\left(0\right)H_{1}\left(0\right)\right)^{2} & \left(H_{3}\left(a\right)H_{1}\left(0\right)\right)^{2} & \cdots & \left(H_{3}\left(3a\right)H_{1}\left(3a\right)\right)^{2} & \left(H_{3}\left(2a\right)H_{1}\left(a\right)\right)^{2}\\
\left(H_{3}\left(0\right)H_{0}\left(0\right)\right)^{2}\\
\left(H_{2}\left(0\right)H_{2}\left(0\right)\right)^{2}\\
\left(H_{2}\left(0\right)H_{1}\left(0\right)\right)^{2} & \vdots &  & \vdots & \vdots\\
\left(H_{2}\left(0\right)H_{0}\left(0\right)\right)^{2}\\
\left(H_{1}\left(0\right)H_{1}\left(0\right)\right)^{2}\\
\left(H_{1}\left(0\right)H_{0}\left(0\right)\right)^{2}\\
\left(H_{0}\left(0\right)H_{0}\left(0\right)\right)^{2} & \left(H_{0}\left(a\right)H_{0}\left(0\right)\right)^{2} & \cdots & \left(H_{0}\left(3a\right)H_{0}\left(3a\right)\right)^{2} & \left(H_{0}\left(2a\right)H_{0}\left(a\right)\right)^{2}
\end{bmatrix}\begin{bmatrix}8 & 0 & \cdots & 0\\
0 & 4 &  & 0\\
\vdots &  & \ddots & \vdots\\
0 & 0 & \cdots & 1
\end{bmatrix}\begin{bmatrix}w_{0,0}\\
w_{1,0}\\
w_{2,0}\\
w_{3,0}\\
w_{4,0}\\
w_{1,1}\\
w_{2,2}\\
w_{3,3}\\
w_{2,1}
\end{bmatrix}=\begin{bmatrix}4!0!\\
3!1!\\
3!0!\\
2!2!\\
2!1!\\
2!0!\\
1!1!\\
1!0!\\
0!0!
\end{bmatrix}.\label{eq:2d orth cond ex}
\end{equation}

Solving the inverse problem for \eqref{eq:2d orth cond ex} gives
us the 2D grid quadrature weights. As with the 1D case, the furthest
node from the centre can be eliminated by choosing the appropriate
grid spacing multiplier for even $N_{q}$. In the above example, $w_{4,0}$
shall be 0. This is one of the reasons that the last node on the last
row has not been shifted to the diagonal. We shall also see that the
weight of this node is equal to that of the furthest node in the 1D
case in the explicit solution. Indeed, solving \eqref{eq:weight elimination condition}
gives us the multiplier required to eliminate that node on the edge.

To derive the modified Smolyak's method, we start from the fundamentals
of Smolyak's method. It is derived by first rewriting the 1D weights
as a telescoping sum as follows
\begin{equation}
\mathbf{w}_{N_{q}}=\mathbf{w}_{0}+\sum_{k=1}^{N_{q}}\mathbf{w}_{k}-\mathbf{w}_{k-1}.
\end{equation}
One can show that summing the tensor products of the weights of arbitrary
dimension $\boldsymbol{W}_{N_{q}-k}$ and differences of 1D weights
of order $k$, mathematically expressed as
\begin{equation}
\boldsymbol{W}_{N_{q}}^{d+1}=\boldsymbol{W}_{N_{q}}^{d}\otimes\mathbf{w}_{0}+\sum_{k=1}^{N_{q}}\boldsymbol{W}_{N_{q-k}}^{d}\otimes\left(\mathbf{w}_{k}-\mathbf{w}_{k-1}\right),\label{eq:smolyak method}
\end{equation}
results in a set of higher rank weights which satisfy the monomial
degree exactness requirement. To see this, first notice that given
the monomial order of the new dimension $n$, the weighted sum of
the monomial from the difference terms for $k>n$ are zero because
the two sets of weights both lead to exact evaluation of the integral.
For $k\le n$, since the sum of the monomial degrees for the rest
of the dimensions is equal to $N_{q}-n$, the weighted sum is the
following telescoping sum,
\begin{align}
\boldsymbol{W}_{N_{q}}^{d+1}\cdot\boldsymbol{\chi}_{N_{q}} & =\left(\boldsymbol{W}_{N_{q-n}}^{d}\cdot\boldsymbol{X}_{N_{q}-n}^{d}\right)\otimes\left(\left(\mathbf{w}_{n}\cdot\mathbf{x}^{n}\right)+\sum_{k=0}^{n-1}\left(\mathbf{w}_{k}-\mathbf{w}_{k}\right)\cdot\mathbf{x}^{n}\right)\label{eq:smolyak telescopic sum proof}\\
 & =\left(\boldsymbol{W}_{N_{q-n}}^{d}\cdot\boldsymbol{X}_{N_{q}-n}^{d}\right)\otimes\left(\mathbf{w}_{n}\cdot\mathbf{x}^{n}\right),
\end{align}
where $\boldsymbol{\chi}$ is the $d+1$ rank monomial tensor of order
$N_{q}$, $\boldsymbol{X}_{N_{q}-n}$ is the $d$ rank monomial tensor
of order $N_{q}-n$, and the weights in the weighted sum of $\boldsymbol{X}_{N_{q}-n}$
are all replaced with the lowest order $N_{q}-n$ because the monomial
order is smaller than or equal to that of the weights. This evaluates
exactly to the moment as all the monomial orders are matched with
the weights.

For the 2D case, setting the weights $\mathbf{w}_{k}$ for both the
differences and the weights on the left side of the tensor product
in \eqref{eq:smolyak method} to 1D weights with $2k+1$ consecutive
nodes as in \eqref{eq:hermite grid quadrature 1D solution} results
in 2D weights positioned in an isosceles triangle on the grid. Attempting
to position the nodes towards the diagonal by setting the Lagrange
interpolation nodes to be at $j=0,\,\left\lfloor \frac{k}{N_{q}-k}\right\rfloor ,\,\left\lfloor \frac{2k}{N_{q}-k}\right\rfloor ,\,\ldots$
for each $\mathbf{w}_{k}$ where $k\le\left\lfloor \frac{N_{q}}{2}\right\rfloor $
in \eqref{eq:smolyak method} leads to denser rows/columns towards
the edge. The trick to eliminate the weights in the in-between nodes
is to modify the Smolyak method to use a linear combination of weights
produced using the in-between nodes for $k\le\left\lfloor \frac{N_{q}}{2}\right\rfloor $.
Rewriting \eqref{eq:smolyak method} to iterate through each designed
term as 
\begin{align}
\mathbf{W}_{N_{q}} & =\begin{cases}
\sum_{k=0}^{\left\lceil \frac{N_{q}}{2}\right\rceil -1}\left(\mathbf{w}_{N_{q}-k}-\mathbf{w}_{N_{q}-k-1}\right)\otimes\overline{\mathbf{w}}_{k}^{N_{q}}+\overline{\mathbf{w}}_{k}^{N_{q}}\otimes\left(\mathbf{w}_{N_{q}-k}-\mathbf{w}_{N_{q}-k-1}\right)+\overline{\mathbf{w}}_{\left\lfloor \frac{N_{q}}{2}\right\rfloor }^{N_{q}}\otimes\overline{\mathbf{w}}_{\left\lfloor \frac{N_{q}}{2}\right\rfloor }^{N_{q}}, & \text{for odd }N_{q},\\
\sum_{k=1}^{\left\lceil \frac{N_{q}}{2}\right\rceil -1}\left(\mathbf{w}_{N_{q}-k}-\mathbf{w}_{N_{q}-k-1}\right)\otimes\overline{\mathbf{w}}_{k}^{N_{q}}+\overline{\mathbf{w}}_{k}^{N_{q}}\otimes\left(\mathbf{w}_{N_{q}-k}-\mathbf{w}_{N_{q}-k-1}\right)+\mathbf{w}_{\left\lfloor \frac{N_{q}}{2}\right\rfloor }\otimes\mathbf{w}_{\left\lfloor \frac{N_{q}}{2}\right\rfloor }, & \text{for even }N_{q},
\end{cases}\label{eq:smolyak for successive elimination}
\end{align}
where $\overline{\mathbf{w}}_{k}^{n}=\sum_{\mathbf{i}\in\Omega_{k,n-k}}b_{k,\mathbf{i}}\mathbf{w}_{k,\mathbf{i}}$,
$\Omega_{k,l}$ denotes the possible combinations of $k$ nodes in
$l$ consecutive nodes, and each $b_{k,\mathbf{i}}\in\mathbb{R}$,
$\sum_{\mathbf{i}\in\Omega_{k,n-k}}b_{k,\mathbf{i}}=1$. The summation
starts from 1 for the even case because from Subsection \ref{subsec:1D-Quadrature-Problem},
we have seen that $\mathbf{w}_{N_{q}}=\mathbf{w}_{N_{q}-1}$ assuming
the conjecture holds. Note that each successive product must contain
a subset of nodes of the previous one, we can design the weights to
eliminate the elements at the nodes we wish to eliminate row by row,
column by column. Because of the symmetry, we consider only the positive
side of the grid. For $k=0$, $\overline{\mathbf{w}}_{0}^{n}$ is
clearly just $\mathbf{w}_{0}=\mathbf{e}_{0}=\begin{bmatrix}1 & 0 & \ldots & 0\end{bmatrix}$.
For the rest, there are $N_{q}-2k$ in-between nodes. Moving the diagonal
node to each of these positions generate exactly $N_{q}-2k$ unknowns
for the matrix system. Because of the symmetry for each $k$ term
in the formulation of \eqref{eq:smolyak for successive elimination},
only the conditions on either the row or the column $N_{q}-k$ need
to be satisfied, there are all together $N_{q}-2k+1$ conditions to
satisfy. Together with the weight for the 1D weights with the nodes
shifted towards the diagonal, and the row for summation of the unknowns,
we have the square matrix of size $N_{q}-2k+1$. This leads to the
system of equations
\begin{equation}
\begin{bmatrix}1 & 1 & \cdots & 1\\
 & \left\{ \bm{\Delta}_{N_{q}-k}\right\} _{N_{q}-k}\left\{ \mathbf{w}_{k,\mathbf{i}_{1}}\right\} _{\left\{ \bar{\mathbf{i}}_{0}\right\} _{0}} &  & 0\\
\left\{ \mathbf{w}_{k}\right\} _{N_{q}-k}\left\{ \bm{\Delta}_{N_{q}-k}\right\} _{\bar{\mathbf{i}}_{0}} &  & \ddots\\
 & 0 &  & \left\{ \bm{\Delta}_{N_{q}-k}\right\} _{N_{q}-k}\left\{ \mathbf{w}_{k,\mathbf{i}_{N_{q}-2k}}\right\} _{\left\{ \bar{\mathbf{i}}_{0}\right\} _{N_{q}-2k-1}}
\end{bmatrix}\begin{bmatrix}b_{k,\mathbf{i}_{0}}\\
b_{k,\mathbf{i}_{1}}\\
\vdots\\
b_{k,\mathbf{i}_{N_{q}-2k}}
\end{bmatrix}=\begin{bmatrix}1\\
\\
-\left\{ \widetilde{\mathbf{W}}_{k}\right\} _{\bar{\mathbf{i}}_{0},N_{q}-k}\\
\\
\end{bmatrix},
\end{equation}
where $\left\{ \mathbf{v}\right\} _{\mathbf{m}}$ denotes a vector
consisting of elements of $\mathbf{v}$ at the indices $\mathbf{m}$,
$\bar{\mathbf{i}}_{0}=\left\{ x:x\notin\mathbf{i}_{0},0<x<N_{q}-k\right\} $,
$\left\{ \mathbf{i}_{0}\right\} _{j}=\left\lfloor \frac{jk}{N_{q}-k}\right\rfloor $,
$\mathbf{i}_{j}=\left(\mathbf{i}_{0}\backslash\left\{ N_{q}-k\right\} \right)\cup\left\{ \bar{\mathbf{i}}_{0}\right\} _{j-1}$,
$\bm{\Delta}_{n}=\mathbf{w}_{n}-\mathbf{w}_{n-1}$, and $\widetilde{\mathbf{W}}_{k}=\sum_{l=0}^{k-1}\bm{\Delta}_{N_{q}-l}\otimes\overline{\mathbf{w}}_{l}^{N_{q}}+\overline{\mathbf{w}}_{l}^{N_{q}}\otimes\bm{\Delta}_{N_{q}-l}$,
except for odd $N_{q}$ with $k=\left\lfloor \frac{N_{q}}{2}\right\rfloor $.
This has the trivial solution
\begin{equation}
\begin{bmatrix}b_{k,\mathbf{i}_{0}}\\
b_{k,\mathbf{i}_{1}}\\
\vdots\\
b_{k,\mathbf{i}_{N_{q}-2k}}
\end{bmatrix}=\begin{bmatrix}c & -cd_{1} & \cdots & -cd_{N_{q}-2k}\\
-cg_{1} & d_{1}+cd_{1}g_{1} & \cdots & cd_{N_{q}-2k}g_{1}\\
\vdots & \vdots & \ddots & \vdots\\
-cg_{N_{q}-2k} & cd_{1}g_{N_{q}-2k} & \cdots & d_{N_{q}-2k}+cd_{N_{q}-2k}g_{N_{q}-2k}
\end{bmatrix}\begin{bmatrix}1\\
\\
-\left\{ \widetilde{\mathbf{W}}_{k}\right\} _{\bar{\mathbf{i}}_{0},N_{q}-k}\\
\\
\end{bmatrix},
\end{equation}
where $d_{j}=\left(\left\{ \bm{\Delta}_{N_{q}-k}\right\} _{N_{q}-k}\left\{ \mathbf{w}_{k,\mathbf{i}_{j}}\right\} _{\left\{ \bar{\mathbf{i}}_{0}\right\} _{j-1}}\right)^{-1}$,
$g_{j}=\frac{\left\{ \mathbf{w}_{k}\right\} _{N_{q}-k}}{\left\{ \bm{\Delta}_{N_{q}-k}\right\} _{N_{q}-k}}\frac{\left\{ \bm{\Delta}_{N_{q}-k}\right\} _{\left\{ \bar{\mathbf{i}}_{0}\right\} _{j-1}}}{\left\{ \mathbf{w}_{k,\mathbf{i}_{j}}\right\} _{\left\{ \bar{\mathbf{i}}_{0}\right\} _{j-1}}}$,
and $c=\left(1-\sum_{j=1}^{N_{q}-2k}g_{j}\right)^{-1}$. Once the
weights $b_{k}$ are solved, one can use them to solve for $k+1$.
And finally, for even $N_{q}$, we add the last contribution $\mathbf{w}_{\left\lfloor \frac{N_{q}}{2}\right\rfloor }\otimes\mathbf{w}_{\left\lfloor \frac{N_{q}}{2}\right\rfloor }$
to the sum.

For odd $N_{q}$, at $k=\left\lfloor \frac{N_{q}}{2}\right\rfloor $,
we can either solve the resulting quadratic equation given by
\begin{multline}
\left(1-b_{0}\right)\left\{ \mathbf{w}_{\left\lceil \frac{N_{q}}{2}\right\rceil }-\mathbf{w}_{\left\lfloor \frac{N_{q}}{2}\right\rfloor ,\mathbf{i}_{0}}b_{0}\right\} _{\left\lfloor \frac{N_{q}}{2}\right\rfloor }\left\{ \mathbf{w}_{\left\lfloor \frac{N_{q}}{2}\right\rfloor ,\mathbf{i}_{1}}\right\} _{\bar{\mathbf{i}}_{0}}+b_{0}\left\{ \mathbf{w}_{\left\lfloor \frac{N_{q}}{2}\right\rfloor ,\mathbf{i}_{0}}\right\} _{\left\lfloor \frac{N_{q}}{2}\right\rfloor }\left\{ \mathbf{w}_{\left\lceil \frac{N_{q}}{2}\right\rceil }-\left(1-b_{0}\right)\mathbf{w}_{\left\lfloor \frac{N_{q}}{2}\right\rfloor ,\mathbf{i}_{1}}\right\} _{\bar{\mathbf{i}}_{0}}\\
+\left(1-b_{0}\right)b_{0}\left\{ \mathbf{w}_{\left\lfloor \frac{N_{q}}{2}\right\rfloor ,\mathbf{i}_{0}}\right\} _{\left\lfloor \frac{N_{q}}{2}\right\rfloor }\left\{ \mathbf{w}_{\left\lfloor \frac{N_{q}}{2}\right\rfloor ,\mathbf{i}_{1}}\right\} _{\bar{\mathbf{i}}_{0}}+\left\{ \widetilde{\mathbf{W}}_{\left\lfloor \frac{N_{q}}{2}\right\rfloor }\right\} _{\bar{\mathbf{i}}_{0},\left\lceil \frac{N_{q}}{2}\right\rceil }=0,
\end{multline}
with $b_{1}=1-b_{0}$, or formulate another linear equation by rewriting
$\left\lfloor \frac{N_{q}}{2}\right\rfloor $-th term of \eqref{eq:smolyak for successive elimination}
in the asymmetric form, by assigning the weights with nodes shifted
towards the diagonal to the left side of the tensor product while
keeping the ones on the right as a linear combination, as
\begin{equation}
\left(\mathbf{w}_{\left\lceil \frac{N_{q}}{2}\right\rceil }-\mathbf{w}_{\left\lfloor \frac{N_{q}}{2}\right\rfloor ,\mathbf{i}_{0}}\right)\otimes\overline{\mathbf{w}}_{\left\lfloor \frac{N_{q}}{2}\right\rfloor }^{N_{q}}+\mathbf{w}_{\left\lfloor \frac{N_{q}}{2}\right\rfloor ,\mathbf{i}_{0}}\otimes\mathbf{w}_{\left\lceil \frac{N_{q}}{2}\right\rceil }.
\end{equation}
Because of the asymmetry, the zero conditions for both the column
and the row need to be satisfied in order to satisfy the condition
for where the non-zero weights are, set out in \eqref{eq:2d orth cond ex}
for instance. While this formulation seems asymmetric in the written
form, it still satisfies the exactness condition given by \eqref{eq:smolyak telescopic sum proof},
and thus the orthogonality condition. Because the solution which satisfies
all the conditions is unique, the final weighted sum shall become
symmetric again. Enforcing the condition for both elements to be zero,
we arrive at the linear system
\begin{multline}
\begin{bmatrix}\left\{ \mathbf{w}_{\left\lceil \frac{N_{q}}{2}\right\rceil }\right\} _{\bar{\mathbf{i}}_{0}}\left\{ \mathbf{w}_{\left\lfloor \frac{N_{q}}{2}\right\rfloor ,\mathbf{i}_{0}}\right\} _{\left\lceil \frac{N_{q}}{2}\right\rceil } & 0\\
0 & \left\{ \mathbf{w}_{\left\lceil \frac{N_{q}}{2}\right\rceil }-\mathbf{w}_{\left\lfloor \frac{N_{q}}{2}\right\rfloor ,\mathbf{i}_{0}}\right\} _{\left\lceil \frac{N_{q}}{2}\right\rceil }\left\{ \mathbf{w}_{\left\lfloor \frac{N_{q}}{2}\right\rfloor ,\mathbf{i}_{1}}\right\} _{\bar{\mathbf{i}}_{0}}
\end{bmatrix}\begin{bmatrix}b_{0}\\
b_{1}
\end{bmatrix}\\
=-\left(\begin{bmatrix}0\\
\left\{ \mathbf{w}_{\left\lfloor \frac{N_{q}}{2}\right\rfloor ,\mathbf{i}_{0}}\right\} _{\left\lceil \frac{N_{q}}{2}\right\rceil }\left\{ \mathbf{w}_{\left\lceil \frac{N_{q}}{2}\right\rceil }\right\} _{\bar{\mathbf{i}}_{0}}
\end{bmatrix}+\begin{bmatrix}\left\{ \widetilde{\mathbf{W}}_{\left\lfloor \frac{N_{q}}{2}\right\rfloor }\right\} _{\bar{\mathbf{i}}_{0},\left\lceil \frac{N_{q}}{2}\right\rceil }\\
\left\{ \widetilde{\mathbf{W}}_{\left\lfloor \frac{N_{q}}{2}\right\rfloor }\right\} _{\bar{\mathbf{i}}_{0},\left\lceil \frac{N_{q}}{2}\right\rceil }
\end{bmatrix}\right),
\end{multline}
where, the last set of unknown weights are given by $\overline{\mathbf{w}}_{\left\lfloor \frac{N_{q}}{2}\right\rfloor }^{N_{q}}=\mathbf{w}_{\left\lfloor \frac{N_{q}}{2}\right\rfloor ,\mathbf{i}_{0}}b_{0}+\mathbf{w}_{\left\lfloor \frac{N_{q}}{2}\right\rfloor ,\mathbf{i}_{1}}b_{1}$.
This again has a trivial solution, which can be used to complete \eqref{eq:smolyak for successive elimination}.

\subsection{Evaluation of Hermite Polynomials}

It is well known that the Hermite polynomials can be evaluated through
the recurrence relationship
\begin{equation}
H_{n+1}=xH_{n}-nH_{n-1}.
\end{equation}
These can be multiplied together through \eqref{eq:hermite tensor}
to form the Hermite polynomial tensor. The higher order difference
method stated by \eqref{lattice boltzman higher order} calls for
the evaluation of the multinomial of Hermite polynomials. This involves
many terms but there is a much more efficient method to evaluate it
because of the isotropy of the Gaussian function. If we show that
the Laplacian Hermite polynomials with $n=0$ and $n=1$ satisfy the
following relationships
\begin{align}
H_{0,m}\left(\mathbf{x}\right) & =rH_{1,m-1}\left(\mathbf{x}\right)-\left(d+2\left(m-1\right)\right)H_{0,m-1}\left(\mathbf{x}\right),\label{eq:hermite laplacian even rec}\\
H_{1,m+1}\left(\mathbf{x}\right) & =rH_{0,m}\left(\mathbf{x}\right)-2mH_{1,m-1}\left(\mathbf{x}\right),\label{eq:hermite laplacian odd recursion}\\
H_{0,0}\left(\mathbf{x}\right) & =1,\\
H_{1,0}\left(\mathbf{x}\right) & =r,\\
\boldsymbol{H}_{1,m}^{\mathbf{e}_{k}}\left(\mathbf{x}\right) & =\frac{\mathbf{x}_{k}}{r}H_{1,m}\left(\mathbf{x}\right),\label{eq:Hermite derivative of Laplacian}
\end{align}
where $r=\left\Vert \mathbf{x}\right\Vert _{2}$, there is no longer
the need to evaluate all the multinomial terms and combine them all
to form the Laplacian. For higher $n$, one can continue applying
the standard recurrence relationship of the Hermite polynomials, starting
from \eqref{eq:Hermite derivative of Laplacian}. In fact, we observe
that the recurrence relationship \eqref{eq:hermite laplacian odd recursion}
is exactly the same as the standard recursion for the Hermite polynomial
of order $2m+1$.

To prove the relationships, first we define the hyperspherical coordinates
for the Cartesian coordinate system of an arbitrary number of dimensions
as
\begin{align}
\mathbf{x}_{0} & =r\prod_{m=1}^{d-1}\sin\theta_{k},\\
\mathbf{x}_{k} & =r\cos\theta_{k}\prod_{\substack{m=k+1}
}^{d-1}\sin\theta_{m},\ \text{for }k=1,\,\ldots,\,d-1.
\end{align}
Clearly, the Hermite polynomials depend only on $r$. Therefore, using
\eqref{eq:hermite laplacian of any order}, one can write in hyperspherical
coordinates, for $n=0,1$,
\begin{align}
H_{0,m}\left(r\right) & =\mathbf{e}^{\frac{r^{2}}{2}}\Delta^{m}\mathrm{e}^{-\frac{r^{2}}{2}},\\
H_{1,m}\left(r\right) & =-\mathrm{e}^{\frac{r^{2}}{2}}\Delta^{m}\frac{\partial}{\partial r}\mathrm{e}^{-\frac{r^{2}}{2}},\label{eq:hermite laplaican 1 definition}
\end{align}
where one can easily show that \eqref{eq:Hermite derivative of Laplacian}
is true by substituting the Jacobian matrix multiplication
\begin{equation}
\begin{bmatrix}\frac{\partial f}{\partial\mathbf{x}_{0}}\\
\frac{\partial f}{\partial\mathbf{x}_{1}}\\
\vdots\\
\frac{\partial f}{\partial\mathbf{x}_{d-1}}
\end{bmatrix}=\begin{bmatrix}\frac{\partial r}{\partial\mathbf{x}_{0}} & \frac{\partial\theta_{1}}{\partial\mathbf{x}_{0}} & \cdots & \frac{\partial\theta_{d-1}}{\partial\mathbf{x}_{0}}\\
\frac{\partial r}{\partial\mathbf{x}_{1}} & \frac{\partial\theta_{1}}{\partial\mathbf{x}_{1}} & \cdots & \frac{\partial\theta_{d-1}}{\partial\mathbf{x}_{1}}\\
\vdots & \vdots & \ddots & \vdots\\
\frac{\partial r}{\partial\mathbf{x}_{d-1}} & \frac{\partial\theta_{1}}{\partial\mathbf{x}_{d-1}} & \cdots & \frac{\partial\theta_{d-1}}{\partial\mathbf{x}_{d-1}}
\end{bmatrix}\begin{bmatrix}\frac{\partial f}{\partial r}\\
\frac{\partial f}{\partial\theta_{1}}\\
\vdots\\
\frac{\partial f}{\partial\theta_{d-1}}
\end{bmatrix}
\end{equation}
into \eqref{eq:hermite laplaican 1 definition}. One can also show
that the Laplacian of an isotropic function $f$ is given by
\begin{equation}
\Delta f=\frac{1}{r^{d-1}}\frac{\partial}{\partial r}r^{d-1}\frac{\partial f}{\partial r},
\end{equation}
using the diagonal metric tensor. Then, for $H_{0,m}$, integrating
by parts gives
\begin{align}
-\frac{1}{r^{d-1}}\frac{\partial}{\partial r}r^{d-1}H_{1,m-1} & =r\mathrm{e}^{\frac{r^{2}}{2}}\Delta^{m-1}\frac{\partial}{\partial r}\mathrm{e}^{-\frac{r^{2}}{2}}+\mathrm{e}^{\frac{r^{2}}{2}}\frac{1}{r^{d-1}}\frac{\partial}{\partial r}r^{d-1}\frac{\partial}{\partial r}\Delta^{m-1}\mathrm{e}^{-\frac{r^{2}}{2}}\\
 & =H_{0,m}-rH_{1,m-1}.
\end{align}
Therefore, it remains to show that $\frac{1}{r^{d-1}}\frac{\partial}{\partial r}r^{d-1}H_{1,m-1}=\left(d+2\left(m-1\right)\right)H_{0,m-1}.$
This is equivalent to the divergence written in the Cartesian coordinate
system as
\begin{align}
\frac{1}{r^{d-1}}\frac{\partial}{\partial r}r^{d-1}H_{1,m-1} & =-\sum_{k=0}^{d-1}\frac{\partial}{\partial\mathbf{x}_{k}}\left(\mathrm{e}^{\frac{r^{2}}{2}}\Delta^{m-1}\frac{\partial}{\partial\mathbf{x}_{k}}\mathrm{e}^{-\frac{r^{2}}{2}}\right),\\
 & =\mathrm{e}^{\frac{r^{2}}{2}}\sum_{k=0}^{d-1}\frac{\partial}{\partial\mathbf{x}_{k}}\Delta^{m-1}\left(\mathbf{x}_{k}\mathrm{e}^{-\frac{r^{2}}{2}}\right)-\mathbf{x}_{k}\Delta^{m-1}\frac{\partial}{\partial\mathbf{x}_{k}}\mathrm{e}^{-\frac{r^{2}}{2}}.\label{eq:hermite laplacian rec proof step2}
\end{align}
Expanding the Laplacian in the second term, we get
\begin{equation}
\sum_{k=0}^{d-1}\frac{\partial}{\partial\mathbf{x}_{k}}\Delta^{m-1}\left(\mathbf{x}_{k}\mathrm{e}^{-\frac{r^{2}}{2}}\right)=\sum_{k=0}^{d-1}\sum_{\sum_{j=0}^{d-1}a_{j}=m-1}\frac{\left(m-1\right)!}{a_{k}!}\frac{\partial^{2a_{k}+1}}{\partial\mathbf{x}_{k}^{2a_{k}+1}}\left(\mathbf{x}_{k}\mathrm{e^{-\frac{\mathbf{x}_{k}^{2}}{2}}}\right)\prod_{\substack{j=0\\
j\neq k
}
}^{d-1}\frac{1}{a_{j}!}\frac{\partial^{2a_{j}}}{\partial\mathbf{x}_{j}^{2a_{j}}}\mathrm{e}^{-\frac{\mathbf{x}_{j}^{2}}{2}}.\label{eq:hermite laplacian rec proof step1}
\end{equation}
Applying the binomial expansion to the product rule, one gets the
identity
\begin{equation}
\frac{\partial^{2a_{k}+1}}{\partial\mathbf{x}_{k}^{2a_{k}+1}}\left(\mathbf{x}_{k}\mathrm{e^{-\frac{\mathbf{x}_{k}^{2}}{2}}}\right)=\mathbf{x}_{k}\frac{\partial^{2a_{k}+1}}{\partial\mathbf{x}_{k}^{2a_{k}+1}}\mathrm{e^{-\frac{\mathbf{x}_{k}^{2}}{2}}}+\left(2a_{k}+1\right)\frac{\partial^{2a_{k}}}{\partial\mathbf{x}_{k}^{2a_{k}}}\mathrm{e^{-\frac{\mathbf{x}_{k}^{2}}{2}}}.
\end{equation}
Substituting it into \eqref{eq:hermite laplacian rec proof step1}
gives
\begin{align}
= & \sum_{\sum_{j=0}^{d-1}a_{j}=m-1}\left(\sum_{k=0}^{d-1}\left(2a_{k}+1\right)\right)\prod_{j=0}^{d-1}\frac{\left(m-1\right)!}{a_{j}!}\frac{\partial^{2a_{j}}}{\partial\mathbf{x}_{j}^{2a_{j}}}\mathrm{e}^{-\frac{r^{2}}{2}}+\sum_{k=0}^{d-1}\mathbf{x}_{k}\Delta^{m-1}\frac{\partial}{\partial\mathbf{x}_{k}}\mathrm{e}^{-\frac{r^{2}}{2}},\\
= & \left(2\left(m-1\right)+d\right)\Delta^{m-1}\mathrm{e}^{-\frac{r^{2}}{2}}+\sum_{k=0}^{d-1}\mathbf{x}_{k}\Delta^{m-1}\frac{\partial}{\partial\mathbf{x}_{k}}\mathrm{e}^{-\frac{r^{2}}{2}},
\end{align}
which can be substituted back into \eqref{eq:hermite laplacian rec proof step2}
to prove \eqref{eq:hermite laplacian even rec}.

For the odd case $H_{1,m}$, using integration by parts, we similarly
have
\begin{align}
\frac{\partial}{\partial r}H_{0,m} & =rH_{0,m}-H_{1,m}.\label{eq:hermite laplacian odd rec 2}
\end{align}
The left hand side, again, can be expressed in the Cartesian coordinate
system as
\begin{align}
\frac{\partial}{\partial r}H_{0,m} & =\frac{r}{\mathbf{x}_{k}}\frac{\partial}{\partial\mathbf{x}_{k}}H_{0,m}\\
 & =\frac{r}{\mathbf{x}_{k}}\sum_{\sum_{j=0}^{d-1}a_{j}=m}\frac{m!}{a_{k}!}\frac{\partial H_{2a_{k}}\left(\mathbf{x}_{k}\right)}{\partial\mathbf{x}_{k}}\prod_{\substack{j=0\\
j\neq k
}
}^{d-1}\frac{H_{2a_{j}}\left(\mathbf{x}_{j}\right)}{a_{j}!}\\
 & =\frac{2mr}{\mathbf{x}_{k}}\sum_{\sum_{j=0}^{d-1}a_{j}=m-1}\frac{\left(m-1\right)!}{a_{k}!}H_{2a_{k}+1}\left(\mathbf{x}_{k}\right)\prod_{\substack{j=0\\
j\neq k
}
}^{d-1}\frac{H_{2a_{j}}\left(\mathbf{x}_{j}\right)}{a_{j}!}\\
 & =\frac{2mr}{\mathbf{x}_{k}}\boldsymbol{H}_{1,m-1}^{\mathbf{e}_{k}}\left(\mathbf{x}\right)=2mH_{1,m-1}\left(\mathbf{x}\right).
\end{align}
Therefore, \eqref{eq:hermite laplacian odd rec 2} is equivalent to
\eqref{eq:hermite laplacian odd recursion}.

\section{\label{sec:Generalisation-to-Fractional}Generalisation to Fractional
Order Laplacian}

The error analysis in Section \ref{subsec:Higher-Order-Approximation}
gives \eqref{eq:detailed error terms}, which states that provided
that $f$ is $\mathrm{C}^{n}$ and the $n+2$ power of Laplacian is
finite at the point of evaluation, one can evaluate its Laplacian
up to $O\left(h^{n}\right)$ accuracy. This result can in fact be
used to generalise to fractional order Laplacian, and further to other
central difference type.

Here, we propose that given a central difference stencil $\boldsymbol{H}$
such that the Laplacian of a function $f$ can be approximated by
the tensor contraction as
\begin{equation}
\Delta f\left(\mathbf{x}\right)=h^{-2}\boldsymbol{H}\cdot\boldsymbol{F}+O\left(h^{2n}\right),\label{eq:general laplacian error}
\end{equation}
where $\boldsymbol{F}$ consists of $f$ evaluated at various symmetric
grid points $\mathbf{x}\pm\mathbf{i}h$, the fractional Laplacian
$-\left(-\Delta^{\frac{\alpha}{2}}\right)f$ can be approximated as
\begin{equation}
-\left(-\Delta^{\frac{\alpha}{2}}\right)f=h^{-\alpha}\boldsymbol{H}_{\alpha}\cdot\boldsymbol{F}+O\left(h^{2n}\right),
\end{equation}
where
\begin{equation}
\boldsymbol{H}_{\alpha}^{\mathbf{m}}=\frac{-1}{\pi^{d}}\dotsint_{0}^{\pi}\prod_{l=0}^{d-1}\cos\left(\mathbf{m}_{l}\mathbf{k}_{l}\right)\left|\sum_{\mathbf{n}\in\boldsymbol{H}^{\mathbf{n}}\neq0}\boldsymbol{H}^{\mathbf{n}}\prod_{l=0}^{d-1}\cos\left(\mathbf{n}_{l}\mathbf{k}_{l}\right)\right|^{\frac{\alpha}{2}}\,\mathrm{d}\mathbf{k}.\label{eq:fractional difference operator}
\end{equation}

The proof can be obtained by realising that \eqref{eq:general laplacian error}
is satisfied if and only if $f$ is at least $\mathrm{C}^{2n}$, and
the discrete-time Fourier transform of the difference operator is
given by
\begin{equation}
\mathcal{F}\left\{ h^{-2}\boldsymbol{H}\right\} \left(\mathbf{k}\right)=h^{-2}\sum_{\mathbf{n}\in\boldsymbol{H}^{\mathbf{n}}\neq0}\boldsymbol{H}^{\mathbf{n}}\prod_{l=0}^{d-1}\cos\left(\mathbf{n}_{l}\mathbf{k}_{l}\right)=-\left|\mathbf{k}\right|^{2}\left(1+\left|\mathbf{k}\right|^{2n}O\left(h^{2n}\right)\right).
\end{equation}
The fractional power of it can be expanded at the origin, using Faadi
Bruno's formula, or its variants \citep{mckiernan1956onthe}, as a
Taylor Series as
\begin{equation}
\left(-\mathcal{F}\left\{ h^{-2}\boldsymbol{H}\right\} \left(\mathbf{k}\right)\right)^{\frac{\alpha}{2}}=\left|\mathbf{k}\right|^{\alpha}\left(1+\left|\mathbf{k}\right|^{2n}O\left(\alpha^{2n}h^{2n}\right)\right).
\end{equation}
The exact coefficients can be evaluated by using FFT, which is faster
than the Faadi Bruno's formula (see \citep{lam2021arbitrary} for
details). Applying inverse DTFT to the above shows that the fractional
operator \eqref{eq:fractional difference operator} is indeed also
$O\left(h^{2n}\right)$. One can also see that since the coefficients
are simply convolution of coefficients of lower order terms, the order
of isotropy remains the same. Therefore, the benefit of improved isotropy
from the generalised lattice Boltzmann method carries over to approximating
the fractional order Laplacian.

The difficulty in such technique is in evaluating the integral given
by \eqref{eq:fractional difference operator}. In \citep{hao2021fractional},
although the fast and efficient method of FFT is applied to evaluate
this integral, we shall see that it will prevent the solution from
converging. Applying FFT is equivalent to applying the trapezoidal
rule. Because of the branch point at the origin, the convergence of
the stencil is limited to $1+\alpha$ for the 1D problem, which when
multiplied by $h^{-\alpha}$ gives an error of $O\left(h\right)$
for $2N+1$ of the stencil coefficients about the origin. This error
has been observed in the numerical experiment in \citep{hao2021fractional}
and seems to be $O\left(h^{2}\right)$ for 2D. However, for the 1D
problem, one can easily show that the second order accurate stencil,
when evaluated with FFT, gives an $O\left(1\right)$ error, meaning
the spectrum of the stencil, and thus the solution will not converge
with respect to the spatial distance used. Using the convolution theorem,
the definition of Dirac comb function, and the closed-form solution
for the stencil \citep{lam2021arbitrary}, one finds that the error
of the stencil evaluated through FFT with $2N+1$ terms is
\begin{equation}
\mathbf{h}_{n}^{\alpha}-\widetilde{\mathbf{h}}_{n}^{\alpha}=-\sum_{\substack{k=-\infty\\
k\neq0
}
}^{\infty}\mathbf{h}_{n+2Nk}^{\alpha}=\sum_{\substack{k=-\infty\\
k\neq0
}
}^{\infty}\left(-1\right)^{n+1}\binom{\alpha}{\frac{\alpha}{2}+n+2Nk}>0.
\end{equation}
Then, the error in the discrete Fourier domain in the origin is
\begin{equation}
\mathrm{FFT}_{2N+1}\left\{ \mathbf{h}^{\alpha}-\widetilde{\mathbf{h}}^{\alpha}\right\} _{0}=\sum_{n=-N}^{N}\left(\mathbf{h}_{n}^{\alpha}-\widetilde{\mathbf{h}}_{n}^{\alpha}\right)>2NhC=O\left(1\right),
\end{equation}
where $C$ is a constant. While this error is not easily generalised
to other combinations of orders and dimensions, we shall see that
this error is similar in other settings through numerical experiments,
and it limits the convergence of the difference method at higher $N$.

One method to avoid the error at the origin is to apply dithering
to the filter. If we can estimate the error amplitude correctly, this
randomises the error and prevents patterns to show up at a specific
band, such as the DC. However, we have not found a way to analyse
the error pattern without significant amount of computation for the
generalisation. Instead, we opt for integration techniques which allow
us to evaluate this integral at higher order convergence rates greater
than the expected convergence rate for solution. With the stencil
error reducing faster than the solution, the convergence of the solution
is no longer hindered.

\subsection{Double Exponential Integration Rule}

One strategy to avoid the branch point is to map it to infinity where
the measure approaches 0, so that the main contribution of the integral
is far away from that branch point. An example of this mapping is
the exponential substitution. With the integrand becoming zero at
the endpoints, the trapezoidal rule converges geometrically. One exponential
rule is the tanh rule, which maps $\left(-\infty,\infty\right)$ to
$\left(-1,1\right)$. However, its convergence is limited by the slow
decay of the measure, requiring a higher number of quadrature points
to cover the range. When linking the the number of quadrature points
with grid distance as in \citep{haber1977thetanh}, it has a convergence
rate of $O\left(\mathrm{e}^{-C\sqrt{N}}\right)$. This is remedied
by the double exponential substitution such as the tanh-sinh rule,
which allows the integral to converge at $O\left(\mathrm{e}^{-CN/\log\left(N\right)}\right)$,
where $C$ is a constant \citep{mori2001thedoubleexponential}. Double
exponential substitution has been applied to solve many integral problems
involving singularity or branch points such as in \citep{beylkin2010approximation},
where the Laplace transform of $t^{\beta}$ is considered.

For Fourier transform, the integrand involves an oscillatory function,
but since it does not introduce any branch cut/point between the lines
$x-\mathrm{i}\frac{\pi}{2}$ and $x+\mathrm{i}\frac{\pi}{2}$, on
which the poles of the tanh-sinh function are located, the convergence
is not affected. However, assuming that we double the number of quadrature
points for every halving of the difference nodal distance, this together
with its exponential convergence still does not answer whether the
evaluated stencil allows the convergence of the solution yet because
we are also doubling the maximum frequency of the oscillatory function
for the additional coefficients. Below, we would like to justify that
given that the $N$-th coefficient of the 1D case of \eqref{eq:fractional difference operator},
given by 
\begin{equation}
\mathbf{h}_{N}^{\alpha}=\frac{-1}{2\pi}\int_{0}^{2\pi}\cos\left(N\omega\right)\left|\sum_{n=-N_{q}}^{N_{q}}\mathbf{h}_{n}\cos\left(n\omega\right)\right|^{\frac{\alpha}{2}}\,\mathrm{d}\omega,
\end{equation}
is approximated by the tanh-sinh rule with $N_{\mathrm{t}}$ nodes
on each side and with error $e$, then the error of the $2N$-th term,
given by
\begin{align}
\mathbf{h}_{2N}^{\alpha} & =-\frac{x^{*}}{N_{t}}\left(\frac{\pi}{4}\left|\sum_{n=-N_{q}}^{N_{q}}\left(-1\right)^{n}\mathbf{h}_{n}\right|^{\frac{\alpha}{2}}+\sum_{k=1}^{N_{t}}\left|\sum_{n=-N_{q}}^{N_{q}}\mathbf{h}_{n}\cos\left(n\tilde{\psi}\left(\frac{x^{*}k}{N_{t}}\right)\right)\right|^{\frac{\alpha}{2}}\cos\left(2N\tilde{\psi}\left(\frac{x^{*}k}{N_{t}}\right)\right)\psi^{\prime}\left(\frac{x^{*}k}{N_{t}}\right)\right),\label{eq:1D de symmetric}
\end{align}
where $\psi\left(x\right)=\tanh\left(\frac{\pi}{2}\sinh\left(x\right)\right)$,
$\tilde{\psi}\left(x\right)=\pi\left(\psi\left(x\right)+1\right)$,
$x^{*}$ is the distance from the origin such that the integrand's
relative magnitude to the origin is smaller than the smallest difference
between two numbers allowed by the floating point precision, is at
max $O\left(e\right)$. The symmetry allows us to only sum from one
side, saving half of the computation. This grid choice as opposed
to relating $N_{t}$ with the grid spacing stems from the fact that
$N_{t}$ is much larger than the typical cases due to the oscillatory
nature of the integrand. The initial guess for $x^{*}$ is based on
the approximation of $\psi^{\prime}$ for large $x$ given by
\begin{equation}
\psi^{\prime}\left(x\right)>\pi\mathrm{e}^{x-\frac{\pi}{2}\mathrm{e}^{x}}.
\end{equation}
The spectrum of the integer order central difference is approximated
by $4\sin^{2}\left(\frac{\omega}{2}\right)$. Then, we have approximately
\begin{equation}
x^{*}\sim\log\left(-\frac{2\left(\log\left(\epsilon\right)-\alpha\log(\pi)\right)}{\pi\left(\alpha+1\right)}\right),\label{eq:limit of xstar}
\end{equation}
where $\epsilon$ is the smallest number which the floating point
number can represent with the exponent being 0. The tolerance $\epsilon$
is further divided by a factor to further account for the underestimation.

To show that the error is, at the maximum, of the same order, we apply
the reproducing kernel
\begin{equation}
K\left(z,w\right)=\frac{\psi^{\prime}\left(z\right)\overline{\psi^{\prime}\left(w\right)}}{1-\psi\left(z\right)\overline{\psi\left(w\right)}},
\end{equation}
which has been employed in \citep{haber1977thetanh} for the error
analysis of the tanh rule, and estimate the $H^{2}$ norm, with respect
to the reproducing kernel, of the error operator given by
\begin{equation}
E^{m}f=\int_{-\infty}^{\infty}\cos\left(m\pi\psi\left(x\right)\right)f\left(x\right)\,\mathrm{d}x-\sum_{k=-\infty}^{\infty}\cos\left(m\pi\psi\left(kh_{\mathrm{e}}\right)\right)f\left(kh_{\mathrm{e}}\right),
\end{equation}
where $f$ is an analytic function on the real line belonging to the
reproducing kernel function space. The norm of a bounded linear operator
in this space is given by \citep{haber1977thetanh}
\begin{equation}
\left\Vert E_{m}\right\Vert ^{2}=E_{\left(z\right)}^{m}E_{\left(w\right)}^{m}K\left(z,w\right).\label{eq:h2 norm of error operator}
\end{equation}
The integral and the sum can be grouped together for analysis by shifting
the integration domain to a line on the imaginary plane as such
\begin{equation}
E^{m}f=\mathfrak{R}\left\{ \int_{\mathrm{i}\beta-\infty}^{\mathrm{i}\beta+\infty}\phi\left(z\right)\cos\left(m\pi\psi\left(z\right)\right)f\left(z\right)\,\mathrm{d}z\right\} ,
\end{equation}
where $0<\beta<\frac{\pi}{2}$, and $\phi\left(z\right)=\left(1-\mathrm{i}\cot\left(\frac{\pi}{h_{\mathrm{e}}}z\right)\right)$,
and $f$ is analytic between the lines $x\pm\mathrm{i}\beta$. For
small $h_{\mathrm{e}}$, $\phi$ can be approximated as
\begin{equation}
\phi\left(x+\mathrm{i}\beta\right)\sim-2\exp\left(\frac{2\pi}{h_{\mathrm{e}}}\left(\mathrm{i}x-\beta\right)\right).
\end{equation}
 The real and imaginary parts of the cosine function are, respectively,
\begin{multline}
\mathfrak{R}\left\{ \cos\left(m\pi\psi\left(x+\mathrm{i}\beta\right)\right)\right\} =\cos\left(\frac{\pi m\sinh\left(\pi\cos\left(\beta\right)\sinh\left(x\right)\right)}{\cos\left(\pi\sin\left(\beta\right)\cosh\left(x\right)\right)+\cosh\left(\pi\cos\left(\beta\right)\sinh\left(x\right)\right)}\right)\cdot\\
\cosh\left(\frac{\pi m\sin\left(\pi\sin\left(\beta\right)\cosh\left(x\right)\right)}{\cos\left(\pi\sin\left(\beta\right)\cosh\left(x\right)\right)+\cosh\left(\pi\cos\left(\beta\right)\sinh\left(x\right)\right)}\right),
\end{multline}
\begin{multline}
\mathfrak{I}\left\{ \cos\left(m\pi\psi\left(x+\mathrm{i}\beta\right)\right)\right\} =-\sin\left(\frac{\pi m\sinh\left(\pi\cos\left(\beta\right)\sinh\left(x\right)\right)}{\cos\left(\pi\sin\left(\beta\right)\cosh\left(x\right)\right)+\cosh\left(\pi\cos\left(\beta\right)\sinh\left(x\right)\right)}\right)\cdot\\
\sinh\left(\frac{\pi m\sin\left(\pi\sin\left(\beta\right)\cosh\left(x\right)\right)}{\cos\left(\pi\sin\left(\beta\right)\cosh\left(x\right)\right)+\cosh\left(\pi\cos\left(\beta\right)\sinh\left(x\right)\right)}\right).
\end{multline}
Now let us consider only the amplitude and ignore the oscillatory
function. The real part can be approximated as
\begin{equation}
\cosh\left(\frac{\pi m\sin(\pi\sin(\beta)\cosh(x))}{\cos\left(\pi\sin\left(\beta\right)\cosh\left(x\right)\right)+\cosh\left(\pi\cos\left(\beta\right)\sinh\left(x\right)\right)}\right)<\frac{1}{2}\left(\mathrm{e}^{m\pi\tan\left(\frac{\pi}{2}\sin\beta\right)\mathrm{e}^{-\frac{\pi}{2}\cos\left(\beta\right)\left(\exp\left(\frac{x^{2}}{2}\right)-1\right)}}+1\right),
\end{equation}
while the imaginary part can be approximated as
\begin{equation}
\sinh\left(\frac{\pi m\sin\left(\pi\sin\left(\beta\right)\cosh\left(x\right)\right)}{\cos\left(\pi\sin\left(\beta\right)\cosh\left(x\right)\right)+\cosh\left(\pi\cos\left(\beta\right)\sinh\left(x\right)\right)}\right)<\mathrm{e}^{m\pi\tan\left(\frac{\pi}{2}\sin\beta\right)\mathrm{e}^{-\frac{\pi}{2}\cos\left(\beta\right)\left(\exp\left(\frac{x^{2}}{2}\right)-1\right)}}-1.
\end{equation}
Since the amplitude of the kernel function behaves like the Gaussian
function, when the real part is multiplied by the kernel, the integral
of the offset part is integrated to a constant, which is $O\left(\mathrm{e}^{-\frac{C}{h_{\mathrm{e}}}}\right)$
when multiplied by $\phi$. For the triple exponential part, which
is also present in the imaginary part, it integrates to $O\left(\mathrm{e}^{Cm}\right)$,
where $0<C<1$, on its own. Therefore, removing the oscillatory functions
from the norm integral \eqref{eq:h2 norm of error operator}, which
is equivalent to setting them to a positive constant function, leading
to an integral of product of positive Gaussian like functions, one
can conclude that, regardless of $\beta$, the integral is $O\left(\mathrm{e}^{C_{1}m-\frac{C_{2}}{h_{\mathrm{e}}}}\right)$,
and so is the norm. Numerical experiment suggests that as $m$ and
$h_{\mathrm{e}}^{-1}$ become very large, $\frac{C_{2}}{C_{1}}$ approaches
1 with $C_{2}>C_{1}$. Therefore, halving of $h$ still allows the
exponential convergence of $2m$-th stencil coefficient. An estimate
of the number of coefficients required by the trapezoidal rule can
be obtained from evaluating the second order stencil, where the exact
solution is known. As the analysis applies to other analytic integrands,
it is also valid for higher order stencils, especially since they
behave similarly.

Another issue that can prevent convergence is the limited numerical
precision. When evaluating the DTFT of the integer order stencil and
the tanh-sinh function near the end point $\pm x^{*}$, the relative
error is much greater if the terms in \eqref{eq:1D de symmetric}
are directly evaluated. Near the end points, the cosine function needs
to be evaluated to a value close to $1$ and be cancelled to obtain
the DTFT value, which is close to zero. However, with the floating
point limited to the number of digits it can represent, the smallest
non-zero DTFT value is limited to that number. If the cosine function
is expressed in terms of the sine function instead, the floating point
can make use of the exponent to represent a much smaller number, allowing
a wider range of DTFT to be evaluated to non-zero values smaller than
$\epsilon$. For tanh-sinh, the series expansion of $\frac{1}{1+x}$
at infinity gives
\begin{equation}
\tanh\left(\frac{\pi}{2}\sinh x\right)=\mathrm{e}^{-\pi\sinh x}-\mathrm{e}^{-2\pi\sinh x}+\mathrm{e}^{-3\pi\sinh x}+\ldots.
\end{equation}
Further expansion allows us to evaluate to even smaller number, but
this is not needed as the solution for $x^{*}$ is around $3.28$
for a constant function and double precision. The attempts to evaluate
these functions more accurately only provide little improvement for
coefficients far from the centre, however, because those coefficients
still require cancellation of the summation terms, which are limited
by the relative error given by $\epsilon$ in \eqref{eq:limit of xstar}.

In the 2D case, it is no longer just a branch point but there are
branch cuts along both the axis, and similarly for 3D, there are branch
surfaces. Applying the tanh-sinh mapping to each of variables will
still map the branching to the infinity, and so the convergence is
not affected. However, we do need the full tensor product, since the
Smolyak algorithm to sparsify the integral operator is not applicable.
This is because each dimension must reach $N$-th term, so we cannot
mix in rules with lower number of terms.

\subsection{Higher order Filon Method}

The tanh-sinh rule clearly has a much higher computational requirement
than FFT. Therefore, we propose that we specify a small number $h_{g}$
away from the origin of the Fourier domain, so that the integral \eqref{eq:fractional difference operator}
is split up into two portions for each dimension, from $0$ to $h_{g}$,
and from $h_{g}$ to $\pi$. The DTFT of the integer order stencil
can now be expanded as a Taylor's series at $h_{g}$, and so a polynomial
based grid quadrature rule can be setup to use FFT coefficients. The
quadrature rule with the sine or cosine function as the weight function
and three point interpolation is termed the Filon-Simpson rule \citep{tuck1967asimple}.
Below, we discuss the generalisation of Filon's approach to higher
order so that the stencil coefficients converge at speeds required.

Evaluating the integral from $0$ to $h_{g}$ means that the trapezoidal
summation is no longer symmetric. While this doubles the amount of
coefficients to be summed, the number of cycles of the cosine functions
are also reduced. This effectively means that the frequency is reduced
when the mapping for the tanh-sinh rule is applied. Following the
deduction from the previous subsection, the number of terms required
to achieve the same error tolerance is also significantly reduced.
Do note that $h_{g}$ should be aligned with the cycles of the cosine
function not only for the previous analysis to directly apply but
also for the FFT coefficients to be applicable. Since the FFT coefficients
are the values of DTFT evaluated at frequencies $\omega=\frac{k\pi}{N}$,
$k=0,\,\ldots,\,2N-1$, this means $h_{g}$ should be $\frac{k\pi}{N}$
for any integer $0\le k\le N$.

For the asymmetrical integrand, since the negative side of the integrand
approaches zero faster than the right hand side, we should find end
points for both sides. For the positive side, the initial guess
\begin{equation}
x_{0}^{\star}=\log\left(-\frac{2}{\pi}\left(\log\left(\epsilon\right)+\alpha\log\left(\frac{\sin\frac{h_{g}}{2}}{\sin h_{g}}\right)\right)\right)
\end{equation}
can be used, where $x^{\star}$ is similarly defined as $x^{*}$ but
the ratio is now between the value of the integrand evaluated at the
points $h_{g}$ and $\frac{h_{g}}{2}$. The initial guess for $x^{*}$
\eqref{eq:limit of xstar} should also be replaced by
\begin{equation}
x_{0}^{*}=\log\left(-\frac{2}{\pi(\alpha+1)}\left(\log\left(\epsilon\right)+\alpha\log\left(\frac{1}{\pi}\sin\frac{h_{g}}{2}\right)\right)\right).
\end{equation}
 Once both $x^{\star}$ and $x^{*}$ are determined, we find a $N_{l}$
such that $\frac{N_{l}x^{\star}}{N_{t}}\ge x^{*}$. And the negative
portion of the summation should terminate at $-N_{l}$.

For the region $h_{g}$ to $\pi$, the integral is further split into
$N-k_{g}$ integration regions, where $k_{g}=\frac{N}{\pi}h_{g}$.
In each region, polynomial interpolation is applied to the spectral
function. For even polynomial order $N_{f}$, 2 consecutive regions
are grouped together so that the interpolation is symmetric over the
range of integration. Then, for each region $m$, $\left(\left(k_{g}+2m\right)h,\left(k_{g}+2m+2\right)h\right)$,
$m=0,\,\ldots,\,\frac{N-k_{g}}{2}$ for even $N_{f}$ and $\left(\left(k_{g}+m\right)h,\left(k_{g}+m+1\right)h\right)$,
$m=0,\,\ldots,\,N-k_{g}$ for odd $N_{f}$, one wishes to find a set
of weights $\mathbf{b}^{n,m}$, for $n=0,\,\ldots\,,N$, such that
\begin{equation}
\int_{\left(k_{g}+qm\right)h}^{\left(k_{g}+q\left(m+1\right)\right)h}\cos\left(nx\right)f\left(x\right)\,\mathrm{d}x=\sum_{j=0}^{N_{f}}f\left(\left(k_{g}+q\left(m+1\right)-1+j-\left\lfloor \frac{N_{f}}{2}\right\rfloor \right)h\right)\mathbf{b}_{j}^{n,m}+O\left(h^{N_{s}}\right),
\end{equation}
where $h=\frac{\pi}{N}$, $q=2$ when $N_{f}$ is even and $q=1$
when $N_{f}$ is odd, $f$ is analytic within the integration limits
including the endpoints, and $N_{s}$ is the convergence order to
be determined. Replacing $f$ with a polynomial of order $N_{f}$,
the summation exactly equals the integral when $\mathbf{b}^{n,m}$
satisfies moments up to degree $N_{f}$, that is,
\begin{equation}
\sum_{j=0}^{N_{f}}\left(\left(k_{g}+q\left(m+1\right)-1+j-\left\lfloor \frac{N_{f}}{2}\right\rfloor \right)h\right)^{k}\mathbf{b}_{j}^{n,m}=\int_{\left(k_{g}+qm\right)h}^{\left(k_{g}+q\left(m+1\right)\right)h}\cos\left(nx\right)x^{k}\,\mathrm{d}x,\quad\text{for }k=0,\,\ldots\,,N_{f}.
\end{equation}
One can shift the integral so that the monomials are 0 at the centre
as such
\begin{equation}
h^{k}\sum_{j=0}^{N_{f}}\left(j-\frac{N_{f}}{2}\right)^{k}\mathbf{b}_{j}^{n,m}=\int_{-\frac{qh}{2}}^{\frac{qh}{2}}\cos\left(n\left(x+c_{m}\right)\right)x^{k}\,\mathrm{d}x,\label{eq:shifted mom for dense filon}
\end{equation}
where $c_{m}=\left(k_{g}+q\left(m+\frac{1}{2}\right)\right)h$. Using
the identity $\cos\left(x\right)=\frac{1}{2}\left(\mathrm{e}^{\mathrm{i}x}+\mathrm{e}^{-\mathrm{i}x}\right)$,
we have
\begin{equation}
\mathbf{J}_{k}^{n,m}=\int_{-\frac{h}{\bar{q}}}^{\frac{h}{\bar{q}}}\cos\left(n\left(x+c_{m}\right)\right)x^{k}\,\mathrm{d}x=\frac{\mathrm{i}^{k-1}}{2n^{k-1}}\left[\Gamma\left(1+k,-\mathrm{i}nx\right)\mathrm{e}^{\mathrm{i}nc_{m}}-\left(-1\right)^{k}\Gamma\left(1+k,\mathrm{i}nx\right)\mathrm{e}^{-\mathrm{i}nc_{m}}\right]_{-\frac{qh}{2}}^{\frac{qh}{2}}.
\end{equation}
The recurrence relationship of the Gamma function allows us to rewrite
the moments in a more computationally friendly form as
\begin{align}
n\mathbf{J}_{k}^{n,m} & =\frac{k}{n}\left(\left(\frac{qh}{2}\right)^{k-1}\cos\left(nb_{m}\right)+\left(-\frac{qh}{2}\right)^{k-1}\cos\left(na_{m}\right)-\left(k-1\right)\mathbf{J}_{k-2}^{n,m}\right)+\left(\frac{qh}{2}\right)^{k}\sin\left(nb_{m}\right)+\left(-\frac{qh}{2}\right)^{k}\sin\left(na_{m}\right),\label{eq:filon moments}
\end{align}
with the initial conditions
\begin{align}
\mathbf{J}_{0}^{n,m} & =\frac{1}{2}\left(\sin\left(nb_{m}\right)-\sin\left(na_{m}\right)\right),\\
\mathbf{J}_{1}^{n,m} & =\frac{\cos\left(nb_{m}\right)-\cos\left(na_{m}\right)}{n^{2}}+\frac{qh}{2}\frac{\sin\left(nb_{m}\right)+\sin\left(na_{m}\right)}{n},
\end{align}
where $a_{m}=c_{m}-\frac{qh}{2}$ , and $b_{m}=c_{m}+\frac{qh}{2}$.
The conditions described by \eqref{eq:shifted mom for dense filon}
can be rewritten in matrix form as
\begin{multline}
\begin{bmatrix}1 &  &  & 0\\
 & h\\
 &  & \ddots\\
0 &  &  & h^{N_{f}}
\end{bmatrix}\begin{bmatrix}1 & 1 & \cdots & 1\\
-\frac{N_{f}}{2} & 1-\frac{N_{f}}{2} & \cdots & \frac{N_{f}}{2}\\
\vdots & \vdots &  & \vdots\\
\left(-\frac{N_{f}}{2}\right)^{N_{f}} & \left(1-\frac{N_{f}}{2}\right)^{N_{f}} & \cdots & \left(\frac{N_{f}}{2}\right)^{N_{f}}
\end{bmatrix}\begin{bmatrix}\mathbf{b}_{0}^{0,0} & \mathbf{b}_{0}^{1,0} & \cdots & \mathbf{b}_{0}^{1,1} & \cdots & \mathbf{b}_{0}^{N,N-k_{g}}\\
\mathbf{b}_{1}^{0,0} & \mathbf{b}_{1}^{1,0} & \cdots & \mathbf{b}_{1}^{1,1} & \cdots & \mathbf{b}_{1}^{N,N-k_{g}}\\
\vdots & \vdots &  & \vdots &  & \vdots\\
\mathbf{b}_{N_{f}}^{0,0} & \mathbf{b}_{N_{f}}^{1,0} & \cdots & \mathbf{b}_{N_{f}}^{1,1} & \cdots & \mathbf{b}_{N_{f}}^{N,N-k_{g}}
\end{bmatrix}\\
=\begin{bmatrix}\mathbf{J}_{0}^{0,0} & \mathbf{J}_{0}^{1,0} & \cdots & \mathbf{J}_{0}^{1,1} & \cdots & \mathbf{J}_{0}^{N,N-k_{g}}\\
\mathbf{J}_{1}^{0,0} & \mathbf{J}_{1}^{1,0} & \cdots & \mathbf{J}_{1}^{1,1} & \cdots & \mathbf{J}_{1}^{N,N-k_{g}}\\
\vdots & \vdots &  & \vdots &  & \vdots\\
\mathbf{J}_{N_{f}}^{0,0} & \mathbf{J}_{N_{f}}^{1,0} & \cdots & \mathbf{J}_{N_{f}}^{1,1} & \cdots & \mathbf{J}_{N_{f}}^{N,N-k_{g}}
\end{bmatrix}.\label{eq:filon system of equations}
\end{multline}
The shifting of the origin of the interpolation polynomial makes the
monomials on the left hand side common for all $m$, so the inverse
only needs to be solved once. The symmetric Vandermonde matrix in
the middle of the left hand side can be solved via Algorithm 1 provided
in \citep{lam2021arbitrary}. 

Next, we analyse the convergence order of the quadrature weights $\mathbf{b}^{n,m}$
with interpolation polynomials of degree $N_{f}$. The interpolation
polynomial can be expressed as 
\begin{equation}
p\left(x\right)=\begin{bmatrix}1 & x-c_{m} & \cdots & \left(x-c_{m}\right)^{N_{f}}\end{bmatrix}\begin{bmatrix}1 &  &  & 0\\
 & h\\
 &  & \ddots\\
0 &  &  & h^{N_{f}}
\end{bmatrix}^{-1}\begin{bmatrix}1 & -\frac{N_{f}}{2} & \cdots & \left(-\frac{N_{f}}{2}\right)^{N_{f}}\\
1 & 1-\frac{N_{f}}{2} & \cdots & \left(1-\frac{N_{f}}{2}\right)^{N_{f}}\\
\vdots & \vdots &  & \vdots\\
1 & \frac{N_{f}}{2} & \cdots & \left(\frac{N_{f}}{2}\right)^{N_{f}}
\end{bmatrix}^{-1}\begin{bmatrix}f\left(c_{m}-\frac{N_{f}}{2}h\right)\\
f\left(c_{m}+\left(1-\frac{N_{f}}{2}\right)h\right)\\
\vdots\\
f\left(c_{m}+\frac{N_{f}}{2}h\right)
\end{bmatrix}.\label{eq:interpolation of f}
\end{equation}
Taylor's series expansion of $f$ about $c_{m}$ gives
\begin{multline}
\begin{bmatrix}f\left(c_{m}-\frac{N_{f}}{2}h\right)\\
f\left(c_{m}+\left(1-\frac{N_{f}}{2}\right)h\right)\\
\vdots\\
f\left(c_{m}+\frac{N_{f}}{2}h\right)
\end{bmatrix}=\begin{bmatrix}1 & -\frac{N_{f}}{2} & \cdots & \left(-\frac{N_{f}}{2}\right)^{N_{f}}\\
1 & 1-\frac{N_{f}}{2} & \cdots & \left(1-\frac{N_{f}}{2}\right)^{N_{f}}\\
\vdots & \vdots &  & \vdots\\
1 & \frac{N_{f}}{2} & \cdots & \left(\frac{N_{f}}{2}\right)^{N_{f}}
\end{bmatrix}\begin{bmatrix}1 &  &  & 0\\
 & h\\
 &  & \ddots\\
0 &  &  & h^{N_{f}}
\end{bmatrix}\begin{bmatrix}\frac{1}{0!} &  &  & 0\\
 & \frac{1}{1!}\\
 &  & \ddots\\
0 &  &  & \frac{1}{N_{f}!}
\end{bmatrix}\begin{bmatrix}f\left(c_{m}\right)\\
f^{\prime}\left(c_{m}\right)\\
\vdots\\
f^{\left(N_{f}\right)}\left(c_{m}\right)
\end{bmatrix}\\
+\begin{bmatrix}\left(-\frac{N_{f}}{2}\right)^{N_{f}+1}\\
\left(1-\frac{N_{f}}{2}\right)^{N_{f}+1}\\
\vdots\\
\left(\frac{N_{f}}{2}\right)^{N_{f}+1}
\end{bmatrix}\frac{h^{N_{f}+1}}{\left(N_{f}+1\right)!}f^{\left(N_{f}+1\right)}\left(c_{m}\right)+O\left(h^{N_{f}+2}\right),
\end{multline}
which can be substituted into \eqref{eq:interpolation of f} to give
\begin{multline}
p\left(x\right)=\begin{bmatrix}1 & x-c_{m} & \cdots & \left(x-c_{m}\right)^{N_{f}}\end{bmatrix}\left(\begin{bmatrix}\frac{1}{0!} &  &  & 0\\
 & \frac{1}{1!}\\
 &  & \ddots\\
0 &  &  & \frac{1}{N_{f}!}
\end{bmatrix}\begin{bmatrix}f\left(c_{m}\right)\\
f^{\prime}\left(c_{m}\right)\\
\vdots\\
f^{\left(N_{f}\right)}\left(c_{m}\right)
\end{bmatrix}+\right.\\
\left.\begin{bmatrix}1 &  &  & 0\\
 & h\\
 &  & \ddots\\
0 &  &  & h^{N_{f}}
\end{bmatrix}^{-1}\begin{bmatrix}1 & -\frac{N_{f}}{2} & \cdots & \left(-\frac{N_{f}}{2}\right)^{N_{f}}\\
1 & 1-\frac{N_{f}}{2} & \cdots & \left(1-\frac{N_{f}}{2}\right)^{N_{f}}\\
\vdots & \vdots &  & \vdots\\
1 & \frac{N_{f}}{2} & \cdots & \left(\frac{N_{f}}{2}\right)^{N_{f}}
\end{bmatrix}^{-1}\begin{bmatrix}\left(-\frac{N_{f}}{2}\right)^{N_{f}+1}\\
\left(1-\frac{N_{f}}{2}\right)^{N_{f}+1}\\
\vdots\\
\left(\frac{N_{f}}{2}\right)^{N_{f}+1}
\end{bmatrix}\frac{h^{N_{f}+1}}{\left(N_{f}+1\right)!}f^{\left(N_{f}+1\right)}\left(c_{m}\right)+O\left(h^{N_{f}+2}\right)\right).
\end{multline}
Comparing to the Taylor's polynomial of $f$ at $c_{m}$, the error
starts at $\left(N_{f}+1\right)$-th term. However, using the symmetry
of each even row and anti-symmetry of each odd row of the Vandemonde
matrix \citep{lam2021arbitrary}, the even degree monomials multiplied
by the $\left(N_{f}+1\right)$-th term are eliminated for even $N_{f}$.
Further expanding the cosine function about $c_{m}$, and integrating
it with the monomials, we find that the local integration error is
\begin{align}
 & h^{N_{f}+1}\int_{-h}^{h}\left(C_{1}x+C_{2}x^{2}+\cdots\right)+O\left(h^{N_{f}+2}\right)\,\mathrm{d}x,\\
= & h^{N_{f}+1}\int_{-h}^{h}\left(C_{2}x^{2}+C_{4}x^{4}+\cdots\right)\,\mathrm{d}x+O\left(h^{N_{f}+3}\right),\\
= & O\left(h^{N_{f}+3}\right),
\end{align}
where $C_{j}$ are bounded constants. Therefore, the global error
order is $N_{s}=N_{f}+q$. Since even order interpolation gives an
additional order of convergence rate, it is superior and preferred.
However, it does impose an additional restriction on $h_{g}$ since
$N-k_{g}$ must now be even. Note that $h_{g}$ cannot be reduced
with increase $N$ or the expansion at $h_{g}$ leads to terms $x^{n}O\left(h_{g}^{\alpha-n}\right)$
which limits convergence to $O\left(h^{1+\alpha}\right)$, which is
the same as the trapezoidal rule. Therefore, the number of terms for
evaluating the tanh-sinh rule must be doubled for each halving of
$h$. For the choice of $N_{f}$, it should satisfy $N_{f}+q>2\max\left(N_{q}-n-N_{c}+1,N_{c}\right)+\alpha$.
For physical problems, $\alpha$ should be less than 2. Therefore,
choosing $N_{f}=2\max\left(N_{q}-n-N_{c}+1,N_{c}\right)$ is sufficient.

Because nodes outside the integration region are used for interpolation,
on the end points, the nodes required may be out of bound of the FFT
coefficients but since the DTFT spectrum is cyclic and symmetric about
$\pi$, one can simply use the mirrored points. The formulation presented
here is not standard compared to the Filon-Simpon method where the
summation uses only nodes within the integration limits. This approach
to integrate a smaller range from the centre of the interpolation
produces a smaller error because of the well-known Runge phenomenon
of Lagrange interpolation. However, this approach does not allow sparsification
of the stencil for higher dimensions as each integration region calls
for neighbouring nodes from the previous region. For 2D integration,
there is not much saving, compared to directly applying the tensor
product on the 1D quadrature weights, as almost half the nodes are
still required, but sparsified stencils may be much more efficient
for 3D. To compute the weights for integrating over the range of end
point nodes, one simply replaces $\frac{qh}{2}$ with $\frac{N_{f}h}{2}$
in \eqref{eq:filon moments} when computing the moments on the right
hand side of \eqref{eq:filon system of equations}. For sparse nodes,
one may either use the Smolyak algorithm described in Subsection \ref{subsec:Multidimensional-Quadrature-Prob}
or solve the system of equations by setting up the moments in a similar
way as \eqref{eq:2d orth cond ex}.

\section{\label{sec:Numerical-Experiment}Numerical Experiment}

In this section, we test the stencil evaluated for 4th order convergence
on two examples and compare the error convergence against stencils
evaluated via FFT. Moreover, we test whether increasing the isotropy
order of the stencil by increasing the quadrature order improves the
isotropy of the error function as predicted. Both examples are isotropic
functions of the radius from the origin. Furthermore, we only consider
the 2D problem. This not only allows us to observe the error pattern
of the stencil, but the analytical solution of their fractional Laplacian
are also easy to find. In the first example, we seek the fractional
Laplacian of
\begin{equation}
f_{1}\left(r\right)=\begin{cases}
\left(1-r^{2}\right)^{\beta} & r\le1,\\
0, & \text{otherwise},
\end{cases}\label{eq:example1}
\end{equation}
where $\beta=6.6$. The solution of its fractional Laplacian in 2D
is given by
\begin{equation}
-\left(-\Delta\right)^{\frac{\alpha}{2}}f_{1}\left(r\right)=\begin{cases}
-2^{\alpha}\frac{\Gamma\left(1+\beta\right)\Gamma\left(1+\frac{\alpha}{2}\right)}{\Gamma\left(1+\beta-\frac{\alpha}{2}\right)}{}_{2}F_{1}\left(1+\frac{\alpha}{2},\frac{\alpha}{2}-\beta;1;r^{2}\right), & r\le1,\\
-\frac{2^{\alpha}}{r^{2+\alpha}}\frac{\Gamma\left(1+\frac{\alpha}{2}\right)}{\left(1+\beta\right)\Gamma\left(-\frac{\alpha}{2}\right)}{}_{2}F_{1}\left(1+\frac{\alpha}{2},1+\frac{\alpha}{2};2+\beta;r^{-2}\right), & r>1,
\end{cases}
\end{equation}
where $_{q}F_{p}$ is the generalised hypergeometric function. Here,
we define the approximation error as
\begin{equation}
E_{i}^{\alpha}=\frac{1}{N_{\mathrm{max}}^{2}}\sum_{j=-N_{\mathrm{max}}}^{N_{\mathrm{max}}}\sum_{k=-N_{\mathrm{max}}}^{N_{\mathrm{max}}}\left|\Delta^{\frac{\alpha}{2}}f\left(r_{j,k}\right)-N_{i}^{-\alpha}\sum_{m=-N_{i}}^{N_{i}}\sum_{n=-N_{i}}^{N_{i}}\mathbf{H}_{n,m}^{\alpha}f\left(r_{j+\frac{N_{\mathrm{max}}}{N_{i}}n,k+\frac{N_{\mathrm{max}}}{N_{i}}m}\right)\right|,
\end{equation}
where $N_{i}=2^{i+4}$, $N_{\mathrm{max}}=2^{8}$, $r_{i,j}=\sqrt{\left(\frac{i}{N_{\mathrm{max}}}\right)^{2}+\left(\frac{j}{N_{\mathrm{max}}}\right)^{2}}$,
and $\mathbf{H}^{\alpha}$ is the 2D central difference stencil. Moreover,
we define the rate of convergence as
\begin{equation}
r_{i}^{\alpha}=\log_{2}\left(\frac{E_{i}^{\alpha}}{E_{i+1}^{\alpha}}\right).
\end{equation}
Tables \ref{tab:Error1} and \ref{tab:Rate1} compare respectively
the approximate error and convergence rate of various stencils. Denoted
as `Sin-FFT' is the stencil defined as
\begin{equation}
\mathbf{H}_{n,m}^{\alpha}=-\frac{1}{4N_{i}^{2}}\sum_{j=0}^{2N_{i}-1}\sum_{k=0}^{2N_{i}-1}\left(4\left(\sin^{2}\frac{\pi j}{N_{i}}+\sin^{2}\frac{\pi k}{N_{i}}\right)\right)^{\frac{\alpha}{2}}\mathrm{e}^{-\mathrm{i}2\pi\left(\frac{nj+mk}{N_{i}}\right)},
\end{equation}
which is the second order stencil described in {[}ref{]}, while `He-FFT'
and `He-Filon' refer to the stencils defined by \eqref{eq:fractional difference operator}
with $N_{c}=2$ and $N_{q}=4$. The integral is respectively approximated
by FFT, and composite double exponential 4th order Filon method. It
can be seen that when the absolute error of solution is small enough
so that the error from the approximation the inverse DTFT by FFT dominates,
neither `Sin-FFT' nor `He-FFT' converges because FFT causes an
$O\left(1\right)$ error as predicted for a function with non-zero
DC in the spectrum (a positive function is guaranteed to have DC).
In fact, a small negative convergence is observed for `Sin-FFT'
possibly due to accumulation of numerical error. When the approximation
error is initially large, for example, in the $\alpha=1.9$ case,
the solution can converge with increasing $N$. `He-FFT' converges
at the same rate as `He-Filon' initially but the rate quickly diminishes.
Similarly for `Sin-FFT', the solution converges at smaller $N$,
but it slows down once the error is saturated. However, `He-Filon'
approaches to the theoretical convergence rate $2N_{c}=4$ as $N$
increases. This example demonstrates that FFT is insufficient for
proper convergence of the solution while higher order integration
method solves this problem.
\begin{table}
\begin{centering}
\begin{tabular}{|c|c|c|c|c|c|c|}
\hline 
\multirow{2}{*}{$i$} & Sin-FFT & He-FFT & He-Filon & Sin-FFT & He-FFT & He-Filon\tabularnewline
\cline{2-7} \cline{3-7} \cline{4-7} \cline{5-7} \cline{6-7} \cline{7-7} 
 & \multicolumn{3}{c|}{$\alpha$=0.1} & \multicolumn{3}{c|}{$\alpha$=0.8}\tabularnewline
\hline 
0 & 0.02381 & 0.01803 & 0.0001463 & 0.01241 & 0.01275 & 0.004889\tabularnewline
\hline 
1 & 0.02381 & 0.01762 & 1.108e-05 & 0.01246 & 0.01251 & 0.0003958\tabularnewline
\hline 
2 & 0.02381 & 0.01717 & 7.343e-07 & 0.01250 & 0.01251 & 2.678e-05\tabularnewline
\hline 
3 & 0.02381 & 0.01670 & 4.666e-08 & 0.01251 & 0.01251 & 1.712e-06\tabularnewline
\hline 
4 & 0.02381 & 0.01618 & 2.929e-09 & 0.01251 & 0.01251 & 1.076e-07\tabularnewline
\hline 
 & \multicolumn{3}{c|}{$\alpha$=1.2} & \multicolumn{3}{c|}{$\alpha$=1.9}\tabularnewline
\hline 
0 & 0.01849 & 0.01996 & 0.01712 & 0.1026 & 0.1227 & 0.1228\tabularnewline
\hline 
1 & 0.008014 & 0.007537 & 0.001432 & 0.02641 & 0.01096 & 0.01099\tabularnewline
\hline 
2 & 0.007499 & 0.007526 & 9.807e-05 & 0.006879 & 0.001041 & 0.0007705\tabularnewline
\hline 
3 & 0.007519 & 0.007525 & 6.289e-06 & 0.001988 & 0.0007552 & 4.971e-05\tabularnewline
\hline 
4 & 0.007523 & 0.007525 & 3.957e-07 & 0.0008402 & 0.0007552 & 3.132e-06\tabularnewline
\hline 
\end{tabular}
\par\end{centering}
\caption{\label{tab:Error1}Error $E_{i}^{\alpha}$ of the approximation of
the fractional Laplacian of \eqref{eq:example1}}

\end{table}
\begin{table}
\begin{centering}
\begin{tabular}{|c|c|c|c|c|c|c|}
\hline 
\multirow{2}{*}{$i$} & Sin-FFT & He-FFT & He-Filon & Sin-FFT & He-FFT & He-Filon\tabularnewline
\cline{2-7} \cline{3-7} \cline{4-7} \cline{5-7} \cline{6-7} \cline{7-7} 
 & \multicolumn{3}{c|}{$\alpha$=0.1} & \multicolumn{3}{c|}{$\alpha$=0.8}\tabularnewline
\hline 
0 & -0.000236 & 0.0336 & 3.72 & -0.00630 & 0.0275 & 3.63\tabularnewline
\hline 
1 & -5.92e-05 & 0.0369 & 3.92 & -0.00391 & 0.000332 & 3.89\tabularnewline
\hline 
2 & -1.48e-05 & 0.0406 & 3.98 & -0.000976 & 6.01e-05 & 3.97\tabularnewline
\hline 
3 & -3.70e-06 & 0.0448 & 3.99 & -0.000244 & 7.45e-05 & 3.99\tabularnewline
\hline 
 & \multicolumn{3}{c|}{$\alpha$=1.2} & \multicolumn{3}{c|}{$\alpha$=1.9}\tabularnewline
\hline 
0 & 1.21 & 1.40 & 3.58 & 1.96 & 3.48 & 3.48\tabularnewline
\hline 
1 & 0.0958 & 0.00213 & 3.87 & 1.94 & 3.40 & 3.83\tabularnewline
\hline 
2 & -0.00375 & 0.000116 & 3.96 & 1.79 & 0.463 & 3.95\tabularnewline
\hline 
3 & -0.000936 & 8.02e-06 & 3.99 & 1.24 & 0.000119 & 3.99\tabularnewline
\hline 
\end{tabular}
\par\end{centering}
\caption{\label{tab:Rate1}Rate of convergence $r_{i}^{\alpha}$ for \eqref{eq:example1}}

\end{table}

In another example, we attempt to approximate the fractional Laplacian
of
\begin{equation}
f_{2}\left(r\right)=\begin{cases}
\left(4r\left(1-r\right)\right)^{n}, & r\le1,\\
0, & \text{otherwise},
\end{cases}\label{eq:example2}
\end{equation}
where $n$ is set to $6$. It has the following solution:
\begin{equation}
-\left(-\Delta\right)^{\frac{\alpha}{2}}f_{2}\left(r\right)=\begin{cases}
\sqrt{\pi}2^{n+\alpha}r^{n-\alpha}\left(\frac{\Gamma\left(\frac{\alpha-1-n}{2}\right)\Gamma\left(2+n\right)r}{\Gamma\left(-\frac{n+1}{2}\right)\Gamma\left(\frac{n}{2}\right)\Gamma\left(\frac{3+n-\alpha}{2}\right)}{}_{4}F_{3}\left(\frac{1-n}{2},1-\frac{n}{2},\frac{3+n}{2},\frac{3+n}{2};\frac{3}{2},\frac{3+n-\alpha}{2},\frac{3+n-\alpha}{2};r^{2}\right)\right.\\
\qquad\left.-\frac{\Gamma\left(\frac{\alpha-n}{2}\right)\Gamma\left(1+n\right)}{\Gamma\left(-\frac{n}{2}\right)\Gamma\left(\frac{1+n}{2}\right)\Gamma\left(\frac{2+n-\alpha}{2}\right)}{}_{4}F_{3}\left(-\frac{n}{2},\frac{1-n}{2},\frac{2+n}{2},\frac{2+n}{2};\frac{1}{2},1+\frac{n-\alpha}{2},1+\frac{n-\alpha}{2};r^{2}\right)\right)\\
\qquad-\frac{2^{1+2n+\alpha}\Gamma\left(\frac{n-\alpha}{2}\right)\Gamma\left(1+n\right)\Gamma\left(1+\frac{\alpha}{2}\right)}{\Gamma\left(1+2n-\alpha\right)\Gamma\left(-\frac{\alpha}{2}\right)}{}_{4}F_{3}\left(\frac{\alpha}{2}-n,\frac{1+\alpha}{2}-n,1+\frac{\alpha}{2},1+\frac{\alpha}{2};1,\frac{1+\alpha-n}{2},1+\frac{\alpha-n}{2};r^{2}\right), & r\le1,\\
-\sqrt{\pi}\frac{2^{\alpha-1}\Gamma\left(1+n\right)\Gamma\left(1+\frac{\alpha}{2}\right)}{r^{2+\alpha}\Gamma\left(\frac{3}{2}+n\right)\Gamma\left(-\frac{\alpha}{2}\right)}{}_{4}F_{3}\left(1+\frac{n}{2},\frac{3+n}{2},1+\frac{\alpha}{2},1+\frac{\alpha}{2};1,\frac{3}{2}+n,2+n;r^{-2}\right), & r>1.
\end{cases}
\end{equation}
Tables \ref{tab:Error2} and \ref{tab:Rate2} list the errors and
convergence rates respectively for the approximation of fractional
Laplacian of \eqref{eq:example2}. Again we observe similar behaviours
from each of the stencils. The consistency displayed here follows
from the previous analysis. 
\begin{table}
\begin{centering}
\begin{tabular}{|c|c|c|c|c|c|c|}
\hline 
\multirow{2}{*}{$i$} & Sin-FFT & He-FFT & He-Filon & Sin-FFT & He-FFT & He-Filon\tabularnewline
\cline{2-7} \cline{3-7} \cline{4-7} \cline{5-7} \cline{6-7} \cline{7-7} 
 & \multicolumn{3}{c|}{$\alpha$=0.1} & \multicolumn{3}{c|}{$\alpha$=0.8}\tabularnewline
\hline 
0 & 0.06178 & 0.04680 & 0.001637 & 0.05065 & 0.08479 & 0.07395\tabularnewline
\hline 
1 & 0.06179 & 0.04574 & 0.0001630 & 0.03263 & 0.03290 & 0.008248\tabularnewline
\hline 
2 & 0.06180 & 0.04458 & 1.200e-05 & 0.03280 & 0.03286 & 0.0006340\tabularnewline
\hline 
3 & 0.06180 & 0.04335 & 7.878e-07 & 0.03284 & 0.03285 & 4.230e-05\tabularnewline
\hline 
4 & 0.06180 & 0.04202 & 4.990e-08 & 0.03285 & 0.03285 & 2.694e-06\tabularnewline
\hline 
 & \multicolumn{3}{c|}{$\alpha$=1.2} & \multicolumn{3}{c|}{$\alpha$=1.9}\tabularnewline
\hline 
0 & 0.1641 & 0.3060 & 0.3022 & 1.484 & 2.815 & 2.816\tabularnewline
\hline 
1 & 0.04761 & 0.04281 & 0.03623 & 0.3957 & 0.3877 & 0.3879\tabularnewline
\hline 
2 & 0.02090 & 0.01996 & 0.002868 & 0.1010 & 0.03283 & 0.03268\tabularnewline
\hline 
3 & 0.01991 & 0.01995 & 0.0001935 & 0.02588 & 0.002914 & 0.002261\tabularnewline
\hline 
4 & 0.01994 & 0.01995 & 1.236e-05 & 0.007044 & 0.002041 & 0.0001456\tabularnewline
\hline 
\end{tabular}
\par\end{centering}
\caption{\label{tab:Error2}Error $E_{i}^{\alpha}$ of the approximation of
the fractional Laplacian of \eqref{eq:example2}}

\end{table}
\begin{table}
\begin{centering}
\begin{tabular}{|c|c|c|c|c|c|c|}
\hline 
\multirow{2}{*}{$i$} & Sin-FFT & He-FFT & He-Filon & Sin-FFT & He-FFT & He-Filon\tabularnewline
\cline{2-7} \cline{3-7} \cline{4-7} \cline{5-7} \cline{6-7} \cline{7-7} 
 & \multicolumn{3}{c|}{$\alpha$=0.1} & \multicolumn{3}{c|}{$\alpha$=0.8}\tabularnewline
\hline 
0 & -0.000280 & 0.0333 & 3.33 & 0.634 & 1.37 & 3.16\tabularnewline
\hline 
1 & -7.00e-05 & 0.0368 & 3.76 & -0.00755 & 0.00185 & 3.70\tabularnewline
\hline 
2 & -1.75e-05 & 0.0405 & 3.93 & -0.00188 & 0.000165 & 3.91\tabularnewline
\hline 
3 & -4.38e-06 & 0.0448 & 3.98 & -0.000469 & 7.97e-05 & 3.97\tabularnewline
\hline 
 & \multicolumn{3}{c|}{$\alpha$=1.2} & \multicolumn{3}{c|}{$\alpha$=1.9}\tabularnewline
\hline 
0 & 1.79 & 2.84 & 3.06 & 1.91 & 2.86 & 2.86\tabularnewline
\hline 
1 & 1.19 & 1.10 & 3.66 & 1.97 & 3.56 & 3.57\tabularnewline
\hline 
2 & 0.0705 & 0.00109 & 3.89 & 1.96 & 3.49 & 3.85\tabularnewline
\hline 
3 & -0.00222 & 6.09e-05 & 3.97 & 1.88 & 0.514 & 3.96\tabularnewline
\hline 
\end{tabular}
\par\end{centering}
\caption{\label{tab:Rate2}Rate of convergence $r_{i}^{\alpha}$ for \eqref{eq:example2}}
\end{table}

\begin{figure}
\begin{centering}
\includegraphics{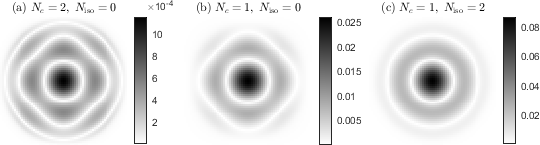}
\par\end{centering}
\caption{\label{fig:ErrorFig1}Error pattern from approximating fractional
Laplacian of \eqref{eq:example1} with $N=64$}
\end{figure}

\begin{figure}
\begin{centering}
\includegraphics{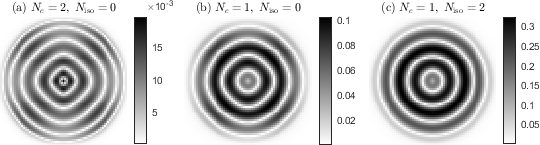}
\par\end{centering}
\caption{\label{fig:ErrorFig2}Error pattern from approximating fractional
Laplacian of \eqref{eq:example2} $N=64$}

\end{figure}

Other than the convergence of solution, another aspect worth looking
into is the isotropy of the stencil. From \eqref{eq:detailed error terms},
the combination of $N_{c}=2$ and $N_{q}=4$ does not lead to isotropic
error terms. Let $N_{\mathrm{iso}}=N_{q}-n-N_{c}+1$, we compare the
stencil evaluated with $N_{\mathrm{iso}}=2$ against those evaluated
with $N_{\mathrm{iso}}=0$ in order to verify the improved isotropy
of the solution. Figures \ref{fig:ErrorFig1} and \ref{fig:ErrorFig2}
illustrate, respectively for $f=f_{1}$ and $f=f_{2}$, the absolute
error defined by
\begin{equation}
\left|\Delta^{\frac{\alpha}{2}}f\left(r_{j,k}\right)-N^{-\alpha}\sum_{m=-N}^{N}\sum_{n=-N}^{N}\mathbf{H}_{n,m}^{\alpha}f\left(r_{j+n,k+m}\right)\right|,
\end{equation}
where $N=64$, $r_{j,k}=\sqrt{\left(\frac{j}{N}\right)^{2}+\left(\frac{k}{N}\right)^{2}}$,
for $j,k=-N,\,\ldots,\,N$. With $N_{\mathrm{iso}}=2$, the highest
order isotropic error term should be 4-th order. This is reflected
in both of the plots, where the error pattern is significantly more
circular than those of stencils with $N_{\mathrm{iso}}=0$. This verifies
the error analysis that the isotropy of the error transfers to the
fractional order stencil. However, this comes at the cost of higher
absolute error when even compared to the stencil with the same error
order. In practice, for diffusion or wave propagation applications,
this isotropy may be preferred over reduced error.


\section{Conclusion}

In conclusion, the differentiation problem of Laplacian has been turned
into multiple integral problems, which requires the application of
multiple quadrature rules, namely the Hermite Gauss quadrature and
quadrature for equidistant nodes, tanh-sinh double exponential substitution
trapezoidal quadrature, and Filon quadrature. By applying the composite
tanh-sinh, and Filon quadrature method, the issue with FFT preventing
the convergence of the solution has been successfully resolved. Moreover,
this method of obtaining higher order stencil is significantly more
efficient than methods presented in \citep{hao2021fractional,lam2021arbitrary}
because the linear combination is applied in the finite space of integer
order stencil. Additionally, the choice of Hermite polynomials allows
us to generate stencils with Gaussian error terms, which are isotropic.
Error isotropy may be more important than absolute error in physical
problems where propagation directions are of great significance.

While we have reviewed in details regarding the lattice Boltzmann
method, including the generalisation to higher order convergence and
exact solutions for 2D stencils, the mystery of the possibility of
the odd order quadrature being able to integrate 1 higher order polynomial
with the appropriate scaling factor for the space variable remains.
Perhaps a future study into the existence of roots of the polynomial
\eqref{eq:weight elimination condition} is worthwhile. Another potential
future study from a mathematical standpoint is a more rigorous analysis
for the convergence of double exponential rule for finite oscillatory
integrals.

\bibliographystyle{elsarticle-num}
\phantomsection\addcontentsline{toc}{section}{\refname}\bibliography{LCD}

\end{document}